\documentclass[leqno,11pt]{article}
\usepackage [dvips]{graphics,psfrag}
\usepackage{amssymb,amsmath}
\usepackage{graphicx}
\usepackage{psfrag}
\usepackage{latexsym}
\usepackage{amsfonts}
\usepackage{url}
\usepackage{rotating}

\newtheorem{theorem}{Theorem}
\newtheorem{proposition}[theorem]{Proposition}
\newtheorem{lemma}[theorem]{Lemma}
\newtheorem{definition}[theorem]{Definition}
\newtheorem{corollary}[theorem]{Corollary}

\newtheorem{example}[theorem]{Example}
\newenvironment{proof}{\noindent{\bf Proof:  }}{\par}
\newcommand{\C}{{\mathbb C}}
\newcommand{\R}{{\mathbb R}}
\newcommand{\Z}{{\mathbb Z}}

\newcommand{\CB}{{\cal B}}

\newcommand{\CH}{{\cal H}}
\newcommand{\CP}{{\cal P}}
\newcommand{\CW}{{\cal W}}
\newcommand{\CL}{{\cal L}}

\renewcommand{\c}{\mathfrak{c}}

\title{Counting Integer Flows in Networks.}

\author{W. Baldoni-Silva, J. A. De Loera, and M. Vergne}
\date{February 17 2003}
\begin{document}
\maketitle

\section*{Introduction}
A {\em network} is a graph with directed edges, with multiple copies
of the edges allowed, and where each node $v$ has an integer value
specified, the so called {\em excess} of $v$, and each arc has an
assigned positive integer value called its {\em capacity}. A {\em
feasible flow} is an assignment of real values to the arcs of the
network so that for any node $v$ the difference between the sum of
values in outgoing arcs minus the sum of values in incoming arcs
equals the prescribed excess of the node $v$ and the capacities of the
arcs are not surpassed. In this paper we study the problem of
effectively counting the number of different {\em integral feasible
flows} in a network. It is well-known that this problem is $\#P$-hard
in the computational category of counting problems \cite{GareyJohnson}
because the problem of counting perfect matchings in bipartite graphs
reduces to it. Despite this bad complexity, concrete applications
abound in graph theory \cite{Jaeger}, representation theory
\cite{kirillov}, and statistics \cite{persi} and thus finding good
methods for attacking concrete examples is of importance. Our goal is
to show that using the algebraic-analytic structure of the problem
allows us to count flows in complicated instances very fast,
surpassing traditional exhaustive enumeration. Continuing the work
started in \cite{baldonivergne} we present effective counting
algorithms from which one can in fact derive counting formulas when
the excess function has parameters.

\begin{figure}[htpb]
\begin{center}
\includegraphics[scale=.6]{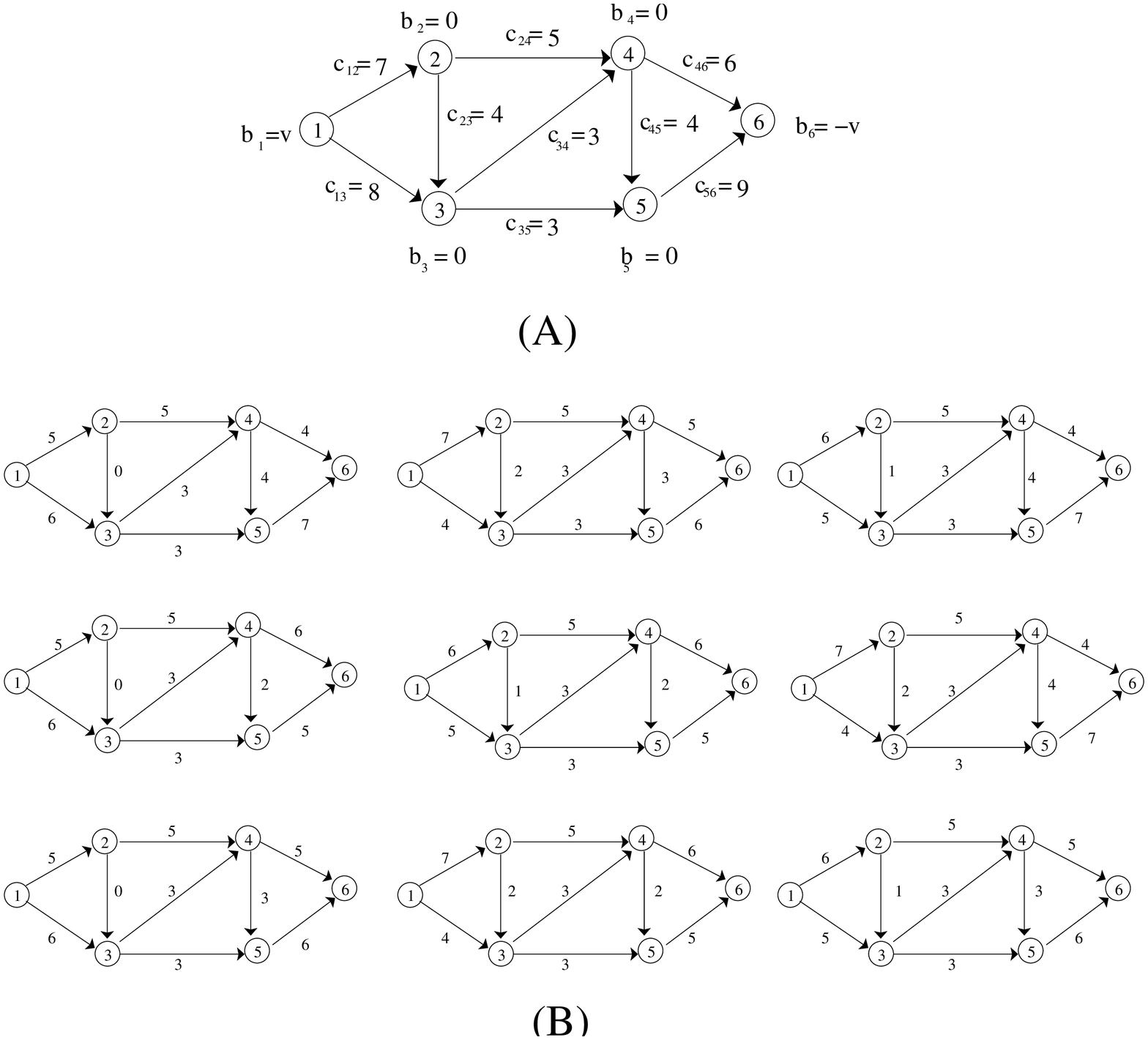}
\caption{counting all maximum flows (part B) of an specific
network (part A)} \label{maxflows}
\end{center}
\end{figure}

The set of all {\it feasible flows} with given excess vector $b$
and capacity vector $c$ is a convex polytope, the well-known {\em
flow polytope}, which is defined by the constraints $\Phi_G
x=b$,\,\,\,$ 0 \leq x \leq c$, where $\Phi_G$ denotes the node-arc
incidence matrix of $G$ (a {\em network matrix}). The incidence
matrix $\Phi_G$ has one column per arc and one row per node. Each
column of $\Phi_G$ has as many entries as nodes. For an arc going
from $i$ to $j$, its corresponding column has zeros everywhere
except at the $i$-th and $j$-th entries. The $j$-th entry, the  
{\em head} of the arrow, receives a $-1$ and the $i$-th entry, 
{\em tail} of the arrow, a $1$. A famous instance
is the {\em max-flow min-cut} problem \cite{schrijver}. This is
the case when $b$ has first entry $v$, last entry $-v$ and $0$
elsewhere. In part (B) of Figure \ref{maxflows} we list all
possible flows with $v=11$, the maximal possible from the network
information specified in part (A).

For us, an important feature of the network incidence matrix $\Phi_G$
is that it is unimodular. We say that the system $\Phi_G$ is {\em
unimodular}, if the columns of $\Phi_G$ span a lattice, denoted by
$\Z\Phi_G$ and, whenever $a$ is in this lattice $\Z\Phi_G$, the
polytope $P(\Phi_G,a)=\{x |\ x\geq 0 :\ \Phi_G x=a\}$ has vertices with
integral coordinates. Even more strongly, network matrices are in fact
totally unimodular matrices \cite{schrijver}, which means that the
lattice generated by their columns is the standard integral lattice
$\Z^n$.  Note that the integral feasible flows are precisely the
integer lattice points inside the flow polytope.

\begin{figure}[htpb]
\begin{center}
\psfrag{1}{$3$}
\psfrag{2}{$1$}
\psfrag{4}{$2$}
\psfrag{a}{$x_1$}
\psfrag{e}{$x_3$}
\psfrag{d}{$x_2$}
\psfrag{c}{$x_4$}
\includegraphics[scale=.7]{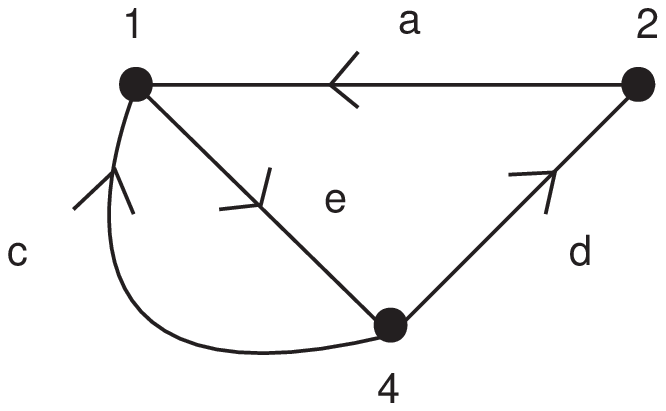}
\caption{Network $G_1$ with nodes $1,2,3$, edges $x_1,x_2,x_3,x_4$, excess
function $b_1=3,b_2=-2,b_3=-1$ and capacity function $c_{x_1}=1 ,
c_{x_2}=1,c_{x_3}=2 ,c_{x_4}=1. $ }
\label{fig:2}
\end{center}

\end{figure}

Here is an example: The node-arc incidence matrix for the graph $G_1$
in Figure \ref{fig:2} is defined by:
$$\Phi_{G_1}=
\begin{pmatrix}
1&-1&0&0\\
0&1&-1&1\\
-1&0&1&-1
\end{pmatrix}.
$$

The equation $\Phi_{G_1}x=b$ reads as the series of equations
$x_1-x_2=b_1$, $x_2-x_3+x_4=b_2$, $-x_1+x_3-x_4=b_3$. These $3$
equations express the fact that, at each node $v\in \{1,2,3\}$,
the difference between the sum of values in outgoing arcs minus
the sum of values in incoming arcs equals the prescribed excess
$b_i$ of the node $v$.  Feasible flows are restricted furthermore
by the conditions $0\leq x_i\leq c_i$.

 The algorithm and formulas deduced here are based on the notion of
{\em total residue} (see Section \ref{keyformulas}), the main concept
involved being the study of rational functions with poles on an
arrangement of hyperplanes.  The enumeration theory we present was
extended to arbitrary rational polyhedra in \cite{Szenesvergne}.  The
particular description we do here is valid for all {\em unimodular
matrices} (again, remember that a matrix $A$ is unimodular if $A$ has
integral coefficients and
the polytope ${\cal P}= \{ x \in {\R}_{+}^m |
Ax= b \},$ has only integral vertices whenever $b$ is in the lattice
spanned by the columns of $A$).

The following lemma implies that it is enough to describe our
counting formulas and techniques for networks without restricted
capacities on the arcs and that have no directed cycles; these are
called {\em acyclic uncapacitated networks}:

\begin{lemma} \label{reducetoacyclic}

 Given a network $G$ with $n$ nodes and $m$ arcs, with capacity $c$
and excess function $b$, there is an acyclic uncapacitated network
$\widehat{G}$ with $n+m$ nodes, $2m$ arcs, and excess function
$\widehat{b}$ (a linear combination of $b$,$c$) such that the
integral flows in both networks are in bijection. The network
$\widehat{G}$ is obtained from $G$ by replacing each arc by two
new arcs as illustrated in the figure below.

 \begin{figure}[h]
 \begin{center}
     \includegraphics[scale=.2]{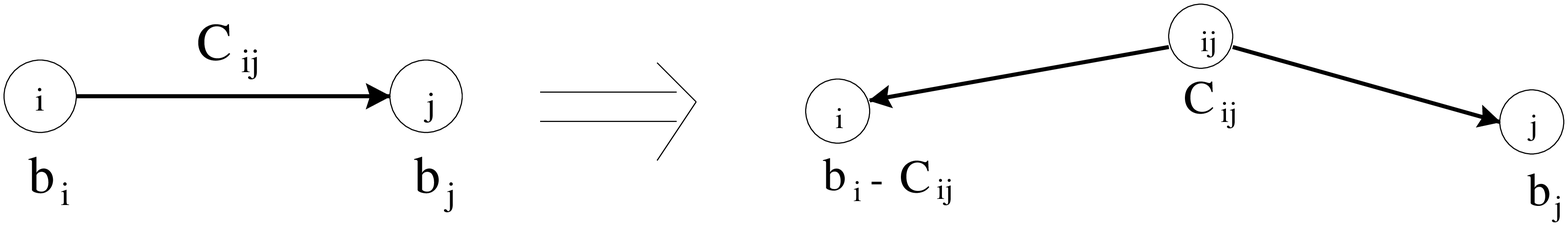}
 \end{center}
 \end{figure}
\label{fig:1}

\end{lemma}

\proof  For the network $G$  with capacity $c,$ the flows are the
solutions of $ \Phi_G  x \, = \, b , \quad 0 \leq x \leq c $ (*).
There is a clear bijection (a projection) between the solutions of
system (*) and the solutions of

$$\left [\begin {array}{cc} {\it \Phi_G}&0\\\noalign{\medskip} I& I \end
{array}\right ]\left [\begin {array}{c} x\\\noalign{\medskip}y
\end {array}\right ] = \left [\begin {array}{c}
b\\\noalign{\medskip}c\end {array}\right ], \quad x, y \, \geq \,
0 .$$

The new enlarged matrix is denoted $\widehat \Phi_G$ and called the
{\em extended network matrix}.  To the network $G$, with its set of
nodes $V$ and its set of arcs $E$, we have associated the new network
$\widehat{G}$ . The set of nodes of $\widehat{G}$ is the disjoint
union of the two sets $V$ and $E$ and the network $\widehat{G}$ is
obtained from $G$ by replacing each arc by two new arcs as illustrated
in the figure above: that is to each $f\in E$ is associated
$f_1=[f,j]$ and $f_2=[f,i]$ where $i$ is the tail of $f$ and $j$ is
the head of $f$. Both arrows $f_1$ and $f_2$ are oriented with their
common tail $\{f\}$ belonging to the set $E$ and their heads $\{i\}$
and $\{j\}$ in the set $V$. Thus $\widehat{G}$ is a directed graph,
with $n+m$ nodes and $2m$ arcs. If $b\in \R^n$ is the excess vector
and $c\in \R^m$ is the capacity vector of the network $G$, we define a
new excess vector $\hat b\in \R^n\oplus \R^m$. The projection of $\hat
b$ on $\R^n$ has coordinates $\hat b_i=b_i-\sum_{f\in E| {
tail}(f)=i}{ capacity}(f)$. The projection of $\hat b$ of $\R^m$ is
the capacity vector $c$. Let $T_G$ be the matrix with one column per
arc and one row per node defined as follows. The column corresponding
to an arc has just {\em one} non zero entry: the {\em tail} of the
arrow receives a $1$.  Then $\Phi_G-T_G$ is the matrix with one column
per arc and just the {\em head} of the arrow receives a $-1$. All
other entries are $0$. Thus
$$\left [\begin {array}{cc} {\it I}&{\it
-T_G}\\\noalign{\medskip} 0& I \end {array}\right ] \left [\begin
{array}{cc} {\it \Phi_G}&0\\\noalign{\medskip} I& I \end
{array}\right ] =\left [\begin {array}{cc} {\it \Phi_G-T_G}&{\it
-T_G}\\\noalign{\medskip} I& I \end {array}\right ] $$ is equal by
definition to the matrix $\Phi_{\widehat G}$, the $m$ first
columns corresponding to new arrows $f_1$, and the last columns
corresponding to new arrows $f_2$. Solutions of
$$\left [\begin {array}{cc} {\it \Phi_G}&0\\\noalign{\medskip} I& I \end
{array}\right ]\left [\begin {array}{c} x\\\noalign{\medskip}y
\end {array}\right ] = \left [\begin {array}{c}
b\\\noalign{\medskip}c\end {array}\right ] $$ are  solutions of
the equation
$$ \Phi_{\widehat G} \left [\begin {array}{c} x\\\noalign{\medskip}y
\end {array}\right ]=\left [\begin {array}{cc} {\it I}&{\it -T_G}\\\noalign{\medskip} 0& I \end
{array}\right ]\left [\begin {array}{c} b\\\noalign{\medskip}c\end
{array}\right ]= \hat b.$$

Thus we obtain a bijection between  feasible flows of the network
$G$ with feasible flows of the uncapacitated network $\widehat G$.
The correspondence assigns to the arc $f_1$ the value $x_f$ and to
the arc $f_2$ the value $y_{f}=c_f-x_f$.

\begin{example}

Consider the network $G_1$ of Figure \ref{fig:2}. Using the
transformation of the previous lemma we would pass from the
capacitated network to the uncapacitated network $G_2$ illustrated
in Figure \ref{fig:4} and the excesses of its nodes are in the
caption.

\begin{figure}[htpb]
        \begin{center}
\psfrag{5}{$5$}
\psfrag{6}{$6$}
\psfrag{7}{$7$}
\psfrag{2}{$2$}
\psfrag{4}{$4$}
\psfrag{1}{$1$}
\psfrag{3}{$3$}
\psfrag{x}{$x_2$}
\psfrag{y}{$(1-x_2)$}
\psfrag{v}{$(1-x_4)$}
\psfrag{k}{$(2-x_3)$}
\psfrag{u}{$x_4$}
\psfrag{w}{$x_1$}
\psfrag{z}{$(1-x_1)$}
\psfrag{s}{$x_3$}

\includegraphics[scale=.6]{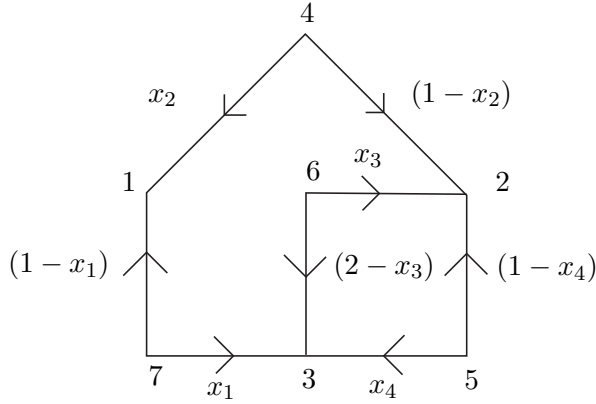}
\caption{Network  $G_2$ with excess $\{2,-4,-3,1,1,2,1\}$ and no
capacities resulting from capacitated network in Figure \ref{fig:2}.}
        \label{fig:4}
\end{center}
\end{figure}

\end{example}

Because of Lemma \ref{reducetoacyclic} and due to interesting
applications in representation theory, it makes sense to focus our
efforts on the special case of uncapacitated  acyclic graphs, and
we do so on Section \ref{kostant}. A particular  case is what
representation theorists would call the {\em Kostant partition
functions} associated to the complete graph $K_n$ with $n$ nodes.
There are many ways to induce an acyclic orientation to the
complete graph, here we take the following convention of
orientation: whenever there is an edge of the graph $G$ between
$i$ and $j$, with $i<j$, then we direct the arrow from $i$ to $j$.

 \begin{example} Consider the complete
graph $G$ on vertices $1,2,3,4.$ In this case,   each vertex is
joined to all the others and the incidence matrix of the network
is

$$\Phi_G=
\begin{pmatrix}   1&1&1&0&0&0\\
{-1}&0&0&1&1&0\\
0&{-1}&0&{-1}&0&1\\
0&0&-1&0&-1&-1
\end{pmatrix}.
$$
\end{example}

Another example of flow polytope is the {\em Pitman-Stanley
polytope} \cite{pitmanstanley} that is constructed starting from a
multiple edge graph:

 \begin{example}\label{ps}
  Consider the graph with vertices $(1,\ldots,n)$ and
edges from $\{i,i+1\}$ and $\{i,n\}$ and the last edge $\{n-1,n\}$
of multiplicity two.  In the case $n=3$ then
$$\Phi_G=
  \begin{pmatrix}
    1 & 1& 0 & 0 \\
    -1 & 0 &1& 1 \\
    0 & -1 & -1 & -1
      \end{pmatrix}.
$$
 \end{example}

  Also, within the class of flow polytopes, we will be investigating
the famous {\em transportation polytopes} \cite{schrijver}. These
polytopes are usually described in terms of $m$ by $n$ real
matrices (denoted here by $M_{m,n}(\R)$): Fix
$c=(c_1,\ldots,c_n)\in \R_+^n$ and $d=(d_1,\ldots,d_m)\in \R_+^m$
such that $\sum_{i=1}^m d_i=\sum_{i=1}^n c_i$ and define
$T_{m,n}(d,c)$ as the set

$$\left\{\begin{array}{llll}
&&x_{ij}\geq 0,1\leq i\leq m,1\leq j\leq n&\\X=\{x_{ij}\} \in
M_{m,n}(\R) &; &\sum_{k} x_{ik}=d_i&1\leq i\leq m\\ &&\sum_{k}
x_{kj}=c_j&1\leq j\leq n\end{array}\right\}.$$

\begin{figure}[h]
 \begin{center}
     \includegraphics[scale=.5]{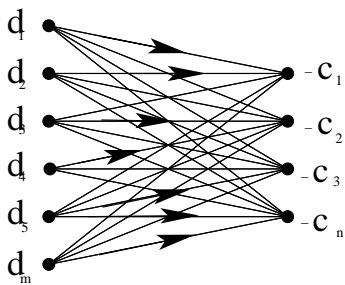}
 \end{center} 
 \caption{The transportation polytopes are network polytopes of
     complete bipartite graphs} \label{bipartite}
 \end{figure}

Then $T_{m,n}(d,c)$ is a polytope called the {\em transportation
polytope} associated for the vectors $d,c$. We can easily see that
this is another flow polytope over a complete bipartite network
$K_{m,n}$ (see Figure \ref{bipartite} where the first $m$ nodes
receive excess values $(d_1,\ldots,d_m)$ and the $n$ nodes in the
second block receive the excess values $(-c_1,-c_2,\ldots,-c_n)$.
The arcs are oriented from the first block to the second. In the
family of transportation polytopes there is a distinguished
member, the {\em Birkhoff polytope} that has been extensively
studied (see for instance the references in the recent paper
\cite{beckpixton}).

It is well-known that the counting formulas of integer flows in a
network come in piecewise polynomial functions (see
\cite{brionvergne,sturmfels}). It is therefore of interest to
understand the regions of validity of each polynomial formula, the
so called {\em chambers}. We dedicate in Section \ref{chambers}
some effort to understand the structure of the chambers and how to
determine the number of chambers. The question of how many chambers
are possible was first raised in \cite{kirillov}. The
combinatorial investigations of the chambers for the partition
functions was initiated by \cite{tatiana1}. See also
\cite{deloeraetal}.

\section{Formulas for the volume and the number of integral
points of flow polytopes.} \label{keyformulas}

In this section, we outline the principles used in the algorithms
we implemented for counting integer flows.  The method is valid
for general convex polytopes \cite{baldonivergne,Szenesvergne},
thus we describe things in a general setting when possible. In
Section \ref{kostant}, we will use particular properties of flow
polytopes associated with  graphs to calculate the counting
formulas.

Let $\Phi$ be an integral  $r$ by $N$ matrix with columns vectors
$ \phi_1,\ldots,\phi_N$. Let $b$ be an $r$-dimensional column
vector and ${\cal{ P}}=\{ x\in {\R}_{+}^N | \Phi x= b \},$ the
rational convex polytope associated to $\Phi$ and $b$. We assume
that $b$ is in the cone $C(\Phi)$ spanned by the non-negative
linear combinations of columns $\phi_1,\phi_2,\ldots,\phi_N$ of
$\Phi$. Without loss of generality we may assume that
$rank(\Phi)=r$. If this is not the case, take the subspace of
$\R^r$ generated by the columns of our matrix and rewrite the
polytope in term of an appropriate rank $k$ matrix of dimension
$k$ by $N$. For example, for the network polytopes the matrices
are not of full rank but deleting one of the rows turns them into
one.

In what follows we assume that $kernel(\Phi) \cap \R^N_+=\{0\}$.
Then $0$ is not in the convex hull of the vectors $\phi_k$ and the
cone $C(\Phi)$ is an acute cone in $\R^r$. For $a\in \R^r$ we
denote by

$P(\Phi,a)=\{(x_1,x_2, \ldots, x_N)\in  \R_+^N | \ \sum_{j=1}^N x_j \phi_j=a\}.$

It is obvious that $P(\Phi,a)$ is a convex polytope determined by
the matrix $\Phi$. Define 
  
 $$v(\Phi,a)=volume(P(\Phi,a)).$$

If $\Phi$ spans a lattice in $\R^r$ and $a$
belongs to this lattice, then define

 $$k(\Phi,a)=\vert P(\Phi,a)\cap \Z^N \vert.$$

Thus $k(\Phi,a)$ is the number of  solutions $(x_1,x_2, \ldots,
x_N)$, in non-negative integers  $x_j$, of the equation
$\sum_{j=1}^N x_j \phi_j={a}.$ The function $k(\Phi,a)$ is called
the {\em vector partition function} associated to $\Phi$. The name
partition comes from the fact that if $\Phi=[e_1,e_1,\ldots ,e_1]$
is the sequence of $N$ times the standard basis vector of $\R$,
then $P_\Phi(a e_1)\cap \Z^N$ is the set of solutions of the
equation $a=x_1+x_2+\cdots+x_N$, that is the partition of the
integer $a$ in $N$ integers. In particular, the function $a \to
k(\Phi,a)$ depends strongly of the multiplicities in the system
$\Phi$. The basic starting observation is

\begin{theorem}
Let $z\in \R^r$ denote a vector in the dual cone to $C(\Phi)$.
Then,

$$ \int_{C(\Phi)} v(\Phi,a) e^{-\langle a,z\rangle } da =\frac{1}{\prod_{\phi \in \Phi} \langle \phi,z\rangle },$$

$$ \sum_{a\in  C(\Phi) \cap \Z^r} k(\Phi,a) e^{-\langle a,z\rangle }=
\frac{1}{\prod_{\phi \in \Phi}1-e^{-\langle \phi,z\rangle }}.$$

\end{theorem}

The goal is to compute the inverses of these two equations. The point
is that one can write efficient formulas for the inversion of Laplace
transforms in terms of residues. In the sequel, we will write
indifferently $\langle \phi,z\rangle $ or $\phi(z)$.

 Let $\Delta^+$ the set $\{\Phi\}$, this means the elements of $\Phi$
 are present {\em without multiplicities}.  We define $\Delta=
 \Delta^+\cup -\Delta^+.$ A subset $\sigma$ of $\Delta^+$ is called a
 {\em basic } subset if $\{\sigma\}$ form a vector space basis of
 $\R^r$.  The {\em chamber complex} is the polyhedral subdivision of
 the cone $C(\Delta^+)$ which is defined as the common refinement of
 the simplicial cones $C(\sigma)$ running over all possible basic
 subsets of $\Delta^+.$ The pieces of this subdivision are called {\em chambers}. We will discuss the chambers in detail, specially how to
 compute the chambers, in Section \ref{chambers}. The important
 fact to remember is that for each chamber there is a quasipolynomial
 formula for $k(\Phi,a)$ and we explain now how to derive the formula
 on a given chamber.

Each $\phi \in \Delta$ determines a linear form on $\C^r$ and a
complex hyperplane $\{z\in \C^r| \phi(z)=0\}$ in $\C^r$. Consider
the hyperplane arrangement $$ \CH_{\C}=\bigcup_{\phi\in
\Delta}\{z\in \C^r | \phi(z)=0\}$$ and let $R_{\Delta}$ denote the
space of rational functions of $z\in \C^r$ with poles on
$\CH_{\C}$. A function in $R_\Delta$ can be written
$P(z)/\prod_{\phi \in \Delta}\phi(z)^{n_\phi}$ where $P$ is a
polynomial function on $r$ complex variables and $n_\phi$ are non
negative integers.  A  subset $\sigma$ of $\Delta$ is called a {\bf
basic subset of} $\Delta$, if the elements $\phi\in \sigma$ form a
vector space basis for $\R^r$. For such $\sigma$, set

$$ f_{\sigma}(z):=\frac{1}{\prod_{\phi\in
\sigma}\phi(z)}. $$

After a linear change of coordinates, the function $f_\sigma$ is
simply $\frac{1}{z_1 z_2\cdots z_r}$ and we denote by $S_{\Delta}$
the subspace of $R_{\Delta}$ spanned by such ``simple" elements
$f_{\sigma}~$. Elements $f_\sigma$ are, in general, {\em not}
linearly independent, as we see in the example below.

\begin{example}

Let $\Delta^+$ be the  set $\Delta^+=\{e_1, e_2, (e_1-e_2)\}.$
Then we have the linear relation
 $$\frac{1}{x y}= \frac{1}{y (x-y)}-
\frac{1}{x (x-y)}$$
 between  elements $f_{\sigma_1}$,$f_{\sigma_2}$,
$f_{\sigma_3}$  with $\sigma_1=\{e_1, e_2\}$, $\sigma_2=\{e_1,
(e_1-e_2)\}$ and  $\sigma_3=\{e_2, (e_1-e_2)\}$ basic subsets  of
$\Delta^+$.
 \end{example}

Partial differentiation $\partial_i$ preserves the space
$R_\Delta$. The key result we need is that there is a well-defined
decomposition of $R_\Delta$ under the action of partial
differentiations, a free module part generated by the basic rational
functions $f_\sigma$,
and a torsion module part, which is unnecessary for calculations
and can be neglected.

\begin{theorem}[Brion-Vergne \cite{brionvergneI}]
The vector space $S_{\Delta}$ is contained in the homogeneous
component of degree $-r$ of $R_{\Delta}$ and we have the direct
sum decomposition

$$R_\Delta=S_\Delta\oplus (\sum_{i=1}^r \partial_i R_\Delta).$$

We call the projection map $$ {  Tres}_\Delta: R_{\Delta}\to
S_{\Delta}
$$ according to this decomposition the {\bf total residue} map.
\end{theorem}
The projection ${  Tres}_\Delta (f)$ of a  function $f$ with poles
on the union of hyperplanes $\CH_{\C}$  depends only of the
smallest hyperplane arrangement $\CH'_{\C}$ containing the poles
of $f$. Therefore we just denote by $ Tres(f)$ the residue of a
rational function $f$ with denominator product of linear forms.
\begin{example}
Observe that if we work in $\R^1$ and $\Delta=\{\pm e_1\}$, then
$R_{\Delta}$ is the space of Laurent series
$$L=\{f(z)=\sum_{k\geq -q}a_k z^k\}.$$ The total residue of a function
$f(z)\in L$ is the function $\frac{a_{-1}}{z}$. The usual residue,
denoted $Res_{z=0} f$, is the constant $a_{-1}$.
\end{example}

 We denote by $\hat R_\Delta$ the obvious extension of $ R_\Delta$,
when we replace the space of polynomial functions on $r$ variables
by the  space of formal power series on $r$ variables.  Let
$F:\C^r \to \C^r$ be an analytic map, such that $F(0)=0$ and
preserving each hyperplane $\phi=0.$ If $f\in \hat R_\Delta$, the
function $(F^*f)(z)=f(F(z))$ is again in $\hat R_\Delta$.  Let ${
Jac}(F)$ be the Jacobian of the map $F$. The function ${ Jac}(F)$
is calculated as follows:  write
$F(z)=(F_1(z_1,z_2,\ldots,z_r),\ldots,F_r(z_1,z_2,\ldots,z_r))$.
Then ${  Jac}(F)=det((\frac{\partial}{\partial z_i}F_j)_{i,j}).$
 We assume ${  Jac}(F)(0)$ does not vanish. For any
$f$ in $\hat R_\Delta$, the following change of variable formula,
which will be useful in our calculations later on, holds in
$S_\Delta$:

$$ Tres (f)={  Tres}({  Jac}(F)(F^*f)).$$

Note that the total residue of a rational function is again a
rational function. By definition, this function can be expressed
as a linear combination of the simple fractions $f_{\sigma}(z)$.
If $f\in S_\Delta$, then ${  Tres}(f)$ is just equal to $f$. We
also know that ${  Tres}$ vanishes on homogeneous rational
functions of degree $m$, whenever $m\neq -r$ and that ${  Tres}$
vanishes on derivatives.  If $f=\frac{P}{\prod_k\langle
\phi_k,z\rangle }$ (with $P$ a polynomial in $r$ variables) has a
denominator product of linear forms $\langle \phi_k,z\rangle $
which do not generate, then it is easy to see that $f$ is a
derivative and the total residue of $f$ is equal to $0$. We are
now ready to fix our notation and recall the key formulas.

\begin{definition}
For $a\in \R^r$, define
 $$J_\Phi(a)(z)=
 {  Tres}(\frac{e^{\langle {a},z\rangle}}{\prod_{k=1}^N
 {\langle \phi_k,z \rangle}})=\frac{1}{(N-r)!}  {  Tres}
 (\frac{\langle  {a},z\rangle^{N-r}}{\prod_{k=1}^N{\langle \phi_k, z\rangle }}
  )$$ and  its ``periodic" version

 $$K_\Phi(a)(z)= {  Tres}( \frac{e^{\langle
{a},z\rangle}}{\prod_{k=1}^N 1-e^{ - \langle\phi_k,z \rangle}}).$$
\end{definition}

The equality:

$${  Tres}(\frac{e^{\langle {a},z\rangle}}{\prod_{k=1}^N
 {\langle \phi_k,z \rangle}})=\frac{1}{(N-r)!}  {  Tres}
 (\frac{\langle  {a},z\rangle^{N-r}}{\prod_{k=1}^N{\langle \phi_k, z\rangle }}
  )$$
follows right away from the fact that the total residue vanishes
on
 homogeneous rational
functions of degree $m$, whenever $m\neq -r$.

By definition, $J_\Phi(a)(z)$ and $K_\Phi(a)(z)$ are rational
functions of $z$ homogeneous in $z$ of degree $-r$. They are
polynomial functions of $a$ of degree $N-r$ and the homogeneous
part
 in $a$ of degree $(N-r)$ in $K_\Phi(a)(z)$ is $J_\Phi(a)(z)$.

\begin{example}

Let us compute $J_\Phi(a)(z)$ and $K_\Phi(a)(z)$ in the case of
the Pitman-Stanley polytope associated to $\Phi_G$ of Example
\ref{ps}. The matrix $\Phi_G$ is a $3$ by $4$ matrix of rank $2$.
Deleting the last row leads to
$$\Phi=
  \begin{pmatrix}
    1 & 1& 0 & 0 \\
    -1 & 0 &1& 1\\
        \end{pmatrix}.
$$

Then  $J_\Phi(a_1,a_2)(z_1,z_2)= {  Tres}\left(\frac{e^{(a_1
z_1+a_2z_2)}}{(z_1-z_2)z_1 z_2^2}\right)$ is
\begin{eqnarray*}
&=&\frac{1}{2!}{  Tres}
 \left(\frac{(a_1 z_1+a_2 z_2)^2}{(z_1-z_2)z_1z_2^2}\right)\\
&=&\frac{a_1^2}{2}{  Tres}
 (\frac{ z_1^2}{(z_1-z_2)z_1z_2^2})+a_1a_2{  Tres}
 (\frac{ z_1z_2}{(z_1-z_2)z_1z_2^2})+\frac{a_2^2}{2}{  Tres}
 (\frac{ z_2^2}{(z_1-z_2)z_1z_2^2})\\
&=&\frac{a_1^2}{2}{  Tres}
 (\frac{ z_1}{(z_1-z_2)z_2^2})+a_1a_2
Tres( \frac{1}{(z_1-z_2)z_2})+\frac{a_2^2}{2}
 Tres(\frac{1}{(z_1-z_2)z_1}).
 \end{eqnarray*}

Now
 $\frac{1}{(z_1-z_2)z_2}$ and $\frac{1}{(z_1-z_2)z_1}$
are simple elements so that they are equal to their respective
total residue. To compute the total residue of $\frac{
z_1}{(z_1-z_2)z_2^2}$,
 we write $z_1$ as  a
linear combination of linear forms in the denominator,
  in order to reduce the degree  of  denominator:
  $$\frac{ z_1}{(z_1-z_2)z_2^2}=
  \frac{(z_1-z_2)+z_2}{(z_1-z_2)z_2^2}=\frac{1}{z_2^2}+\frac{1}{(z_1-z_2)z_2}.$$
The total residue of $\frac{1}{z_2^2}$ is $0$, as
$\frac{1}{z_2^2}= -\frac{\partial}{\partial_{z_2}}\frac{1}{z_2}$
is a derivative, thus $Tres(\frac{
z_1}{(z_1-z_2)z_2^2})=\frac{1}{(z_1-z_2)z_2}$. We finally obtain:
$$J_\Phi(a_1,a_2)(z_1,z_2)=\frac{1}{2}
\frac {a_1^2+2 a_{1}a_{2}} {(z_{1} - z_{2})z_{2}}+ \frac {1}{2}
\frac {a_{2}^{2}}{(z_{1} - z_{2}) \,z_{1}}.$$

We now compute:

$$K_\Phi(a_1,a_2)(z_1,z_2)=
{  Tres}\left(\frac{e^{(a_1
z_1+a_2z_2)}}{(1-e^{-(z_1-z_2)})(1-e^{-z_1}) (1-e^{-
z_2})^2}\right).$$ This is
$${  Tres}\left(\frac{1}{(z_1-z_2)z_1z_2^2}
e^{(a_1 z_1+a_2z_2)}
\frac{(z_1-z_2)}{(1-e^{-(z_1-z_2)})}\frac{z_1}
{(1-e^{-z_1})}\frac{z_2^2}{(1-e^{- z_2})^2}\right).$$

We replace the analytic function
$$e^{(a_1 z_1+a_2z_2)} \frac{(z_1-z_2)}{(1-e^{-(z_1-z_2)})}\frac{z_1}
{(1-e^{-z_1})}\frac{z_2^2}{(1-e^{- z_2})^2}$$ by its Taylor series
at $z_1=0,z_2=0$, and keep only its  term $N(a_1,a_2)(z_1,z_2)$ of
homogeneous degree $2$ in $z_1,z_2$ which is

$$({\displaystyle \frac {5}{12}}  + {a_{1}} +
{\displaystyle \frac {1}{2}} \,{a_{1}}^{2})\,{z_{1}}^{2} + (
{\displaystyle \frac {7}{12}}  + {a_{2}} + {\displaystyle \frac {
1}{2}} \,{a_{1}} + {a_{1}}\,{a_{2}})\,{z_{1}}\,{z_{2}} + (
{\displaystyle \frac {1}{2}} \,{a_{2}} + {\displaystyle \frac {1
}{2}} \,{a_{2}}^{2})\,{z_{2}}^{2}.$$

Thus $K_\Phi(a_1,a_2)(z_1,z_2)$ is equal to
$$Tres\left( \frac{N(a_1,a_2)(z_1,z_2)}{(z_1-z_2)z_1z_2^2}\right).$$

Arguing as for $J_\Phi$, we finally obtain
 that $K_\Phi(a_1,a_2)(z_1,z_2)$ is equal to

$$
\frac{1}{2} \frac {a_1^2+2 a_{1}a_{2}+3 a_1+2 a_2+2} {(z_{1} -
z_{2})z_{2}}+ \frac {1}{2} \frac {a_{2}^{2}+a_2}{(z_{1} - z_{2})
\,z_{1}}.$$

\end{example}

We are now ready to write the formulas to compute the volume and
number of integral points. See \cite[Section 2]{baldonivergne} for
details. To each chamber $\c$ of the subdivision of $C(\Delta^+)$
is associated a linear form $f\to \langle \langle
\c,f\rangle\rangle$ on $S_\Delta$. If the system $\Phi$ is
unimodular, as is the case for networks, it takes value $1$ or $0$
on  $f_\sigma$ whether or not $\c$ is contained in $C(\sigma)$.

\begin{theorem}[Baldoni-Vergne \cite{baldonivergne}]\label{formule}

Let $\c$ be a chamber of the subdivision of $C(\Delta^+)$

\begin{enumerate}
\item For $ a\in \overline{\c}$, the volume of $P(\Phi,a)$ is given by
$$v(\Phi,a)=\langle  \langle  \c, J_\Phi(a)\rangle\rangle.$$

\item If the system $\Phi$ is
unimodular, then for $ a\in \overline{\c} \cap \Z\Phi$, the number
of
 integral points in $P(\Phi,a)$ is given by

$$k(\Phi,a)= \langle\langle \c, K_\Phi(a)\rangle\rangle.$$

\item The function $a\mapsto v(\Phi,a)$ is polynomial on a
chamber $\c$.

\item If the system $\Phi$ is unimodular, as is the case for networks,
 the Ehrhart function $a\mapsto k(\Phi,a)$ is polynomial on a
specified neighborhood of a chamber $\c$.
\end{enumerate}

\end{theorem}

A more general formula for arbitrary $\Phi$ spanning a lattice
$\Z\Phi$ in $\R^r$ is given in \cite{Szenesvergne}. Now, the
question is how to apply these two formulas for the computations
with flow polytopes. The calculation of total residues will
simplify considerably.

\section{Counting Integer Flows in Networks}
\label{kostant}

In this section we will focus on flow polytopes for acyclically
directed graphs. We already justified in the introduction this
makes sense, as other networks can be reduced to acyclic
uncapacitated networks. Consider a $r+1$ real dimensional vector
space.  Let $A_r^+$ (the positive root system of $A_r$) be defined by

$$A_r^+=\{(e_i-e_j)| 1\leq i< j\leq (r+1)\}.$$

Consider $E_r$   the vector space spanned by the elements
$(e_i-e_j)$, then
 $$E_r=\{ a\in \R^{r+1}|
 a=a_1 e_1+\cdots+a_{r} e_r+a_{r+1} e_{r+1} \,\,\,\text{with}\, \,
a_1+a_2+\cdots+a_r+a_{r+1}=0\}.$$
 The vector space $E_r$ is of dimension $r$ and the map
\begin{equation}\label{iden}
f: \ \R^r \longrightarrow E_r
\end{equation}

defined by
$$a=(a_1,a_2,\ldots, a_r) \longmapsto {\bf a} =a_1 e_1+\cdots+a_{r}
e_r-(a_1+\cdots + a_r)e_{r+1}$$ explicitly provides an isomorphism
of $E_r$ with the Euclidean space $\R^r.$  Let, as before,
$\Phi=\{\phi_1, \phi_2,\ldots,\phi_N\}$ denote a sequence of
non-zero linear forms belonging to $A_r^+$. We assume that the
vector space spanned by $\Phi$ is $E_r$. This sequence  is
completely specified by the multiplicity  $m_{i,j}$  of the vector
$e_i-e_j$ in $\Phi$. Explicitly for the transportation polytope
$T_{m,n}(d,c)$, if we denote by $\Phi_{m,n}\subset A^+_{m+n-1}$
the roots  associated to it, then we have $\Phi_{m,n}=\{(e_i-e_j) |
1\leq i<m,\ m+1\leq j\leq m+n\}$  and  thus $m_{i,j}=1$ if $1\leq
i\leq m, \ m+1 \leq j \leq m+n$, $m_{i,j}=0$ otherwise.

It is clear that the polytope $P({\Phi},a)$ is the polytope associated
to the uncapacitated network with $(r+1)$ nodes, where the arc
$i\mapsto j$ ($i<j$) appears $m_{i,j}$ times ($m_{i,j}$ can be $0$ for
some arcs), and with excess function $a_i$ at each node $1,2,\ldots,
r$ and $-(a_1+a_2+\cdots+a_r)$ at the last node $r+1$.  Indeed we have
seen in Remark 3 that the columns of the matrix corresponding to
$P({\Phi},a)$ are vectors of the form $e_i-e_j$ for some $i$ and $j.$

The hyperplane arrangement  (setting $z_{r+1}=0$) generated by
$A_r^+$ is given by the following set of hyperplanes:

 $$\{ z_i | 1\leq i\leq r\} \cup \{(z_i-z_j)| 1\leq i< j\leq r\}.$$

A function in $R_{A_r}$ is  thus a rational function
$f(z_1,z_2,\ldots,z_r)$ on $\C^r$, with poles on the hyperplanes
$z_i=z_j$ or $z_i=0$. The following result is proved by induction
in \cite{baldonivergne}, Proposition 14.

\begin{lemma}

 Let $\Sigma_r$  be the set of permutations on $\{1,2,\ldots,r\}$ and $f_{\pi}$, $f_w, w\in  \Sigma_r$ be defined
by
$$f_{\pi}(z_1,z_2,\ldots,z_r)=\frac{1}{(z_1-z_2)(z_2-z_3)\cdots
(z_{r-1}-z_{r})z_r} $$ and  $$f_w(z_1,\ldots,z_r)=w\cdot
f_{\pi}(z_1,\ldots,z_r)=\frac{1}{(z_{w(1)}-z_{w(2)})(z_{w(2)}-z_{w(3)})\cdots
(z_{w(r-1)}-z_{w(r)})z_{w(r)}}$$ then

\begin{eqnarray}\label{dim}
{  dim}\  S_{A_r}=r!\qquad \qquad \qquad\qquad \qquad \qquad \qquad \qquad \qquad \qquad \qquad  \ \ \   \\
\mbox{and} \qquad \qquad \qquad\qquad \qquad \qquad \qquad \qquad
\qquad \qquad \qquad\qquad \qquad\qquad\qquad \nonumber
\\\{f_w(z_1,\ldots,z_r)=w\cdot f_{\pi}(z_1,\ldots,z_r),  \ \ w\in \Sigma_r\} \ \ \qquad\qquad
\qquad\qquad\\
\mbox{is a basis for}\ S_{A_r} \ \ \qquad \qquad \qquad
\qquad\qquad \qquad\qquad \qquad \qquad \qquad \qquad\qquad
\nonumber
\end{eqnarray}

\end{lemma}

   The cone  $C(A_r^+)$  generated by positive
roots is the cone $a_1\geq 0$, $a_1+a_2\geq 0$,\ldots$,
a_1+a_2+\cdots +a_r\geq 0$.  We denote by $\c^+$ the  open set  of
$C(A_r^+)$ defined by
 $$\c^+=\{a\in C({A_r}^+)\,\, | \,\,a_i>0 ,\,i=1,\ldots,r\}.$$ It is a
chamber of our subdivision, and will be called the {\em nice
chamber}. The importance of this chamber is that its ``permutations''
form a "basis" for the formulas that express volume and number of
integral points. If $\c$ is a chamber for $C(\Phi)$ then there exists
a unique chamber of $C(A_r^+)$ that contains $\c.$

\begin{definition}[\cite{baldonivergne}]
Let $m_{i,j} $ ($i<j$) be the multiplicity of the vector $e_i-e_j$
in $\Phi$ (i.e. this is the number of times the arc $i,j$ is
present in the network). Let $N=\sum_{i,j}m_{i,j}$ the total
number of arcs. We explicitly write down the functions $J_\Phi(a)$
and $K_\Phi(a)$ for our choice of $\Phi$, $a$. Recalling that
$z_{r+1}=0$, we have that

\begin{itemize}
\item
$J_\Phi(a)(z_1,\ldots,z_r)=\frac{1}{(N-r)!}{
Tres}\left(\frac{(a_1z_1+\cdots+a_rz_r)^{N-r}}
{z_1^{m_{1,r+1}}z_2^{m_{2,r+1}}\cdots z_r^{m_{r,r+1}} \prod_{1\leq
i< j\leq r}(z_i-z_j)^{m_{ij}}}\right),$
\item $K_\Phi(a)(z_1,\ldots,z_r)={  Tres}\left(\frac{e^{a_1  z_1}e^{a_2 z_2}\cdots e^{a_r z_r}}
{\prod_{i=1}^r(1-e^{-z_i})^{m_{i,r+1}}\prod_{1\leq i< j\leq
r}(1-e^{-(z_i-z_j)})^{m_{ij}}}\right).$
\end{itemize}

\end{definition}

We now write these functions in two specific examples.

\begin{example}
We consider the polytope associated to a complete bipartite graph
with $3$ nodes on each side. Recall that in this case the matrix
that determines the polytope is given by the vectors
$\Phi=\{e_1-e_4,e_1-e_5,e_1-e_6,e_2-e_4,e_2-e_5,e_2-e_6,e_3-e_4,e_3-e_5,e_3-e_6\}.$
So

\[ \left\{\begin{array}{lll}
m_{i,j}=1 &\qquad\mbox{if}&   1\leq i\leq 3 \  \mbox{and} \ 4\leq j\leq 6\\
0&&\qquad\mbox{otherwise}\qquad\end{array}\right.
\]

and
\begin{itemize}
\item $J_\Phi(a)(z_1,\ldots,z_5)=\frac{1}{4!}{  Tres}\left(\frac{(a_1z_1+a_2z_2+a_3z_3+a_4z_4+a_5z_5)^4}
{z_1 z_2 z_3\prod_{ 1\leq i\leq 3\atop 4\leq j\leq 5 }
(z_i-z_j)}\right),$
\item $K_\Phi(a)(z_1,\ldots,z_5)={  Tres}\left(\frac{e^{a_1  z_1}e^{a_2 z_2}e^{a_3 z_3}e^{a_4 z_4}e^{a_5
z_5}} {\prod_{i=1}^3(1-e^{-z_i})\prod_{ 1\leq i\leq 3\atop 4\leq
j\leq 5 }(1-e^{-(z_i-z_j)})}\right).$
\end{itemize}
\end{example}

\begin{example}
We consider the polytope determined by the complete graph  $K_5$,
in other words $\Phi=A^+_4.$ We obtain

\begin{itemize}
\item $J_\Phi(a)(z_1,\ldots,z_4)=\frac{1}{6!}{  Tres}\left(\frac{(a_1z_1+a_2z_2+a_3z_3+a_4z_4)^6}
{z_1 z_2 z_3z_4\prod_{ 1\leq i< j\leq 4 } (z_i-z_j)}\right),$
\item $K_\Phi(a)(z_1,\ldots,z_4)={  Tres}\left(\frac{e^{a_1  z_1}e^{a_2 z_2}e^{a_3 z_3}e^{a_4 z_4}}
{\prod_{i=1}^4(1-e^{-z_i})\prod_{1\leq i< j\leq
4}(1-e^{-(z_i-z_j)})}\right).$
\end{itemize}
\end{example}

In handling the formulas that we have for computing the volume
and the number of integral points, the first problem is that of
computing the total residue. This is in general a very difficult
task. On the other hand, as we have seen, there is a very nice
basis in $S_{A_r}$ and this will allow us to rewrite the formulas
in terms of iterated residue, which are certainly more tractable.
The point is that one needs to find some, but not all, simplicial
cones that contain the chamber determined by $a$. This is a step
that allows the complexity of the algorithm to be reduce. We are
now going to introduce the iterated residue for $A_r.$

 Recall that, via the identification (\ref{iden}) of $E_r$ with $\R^r$,
a function in $R_{A_r}$ is a rational function
$f(z_1,z_2,\ldots,z_r)$ on $\C^r$, with poles on the hyperplanes
$z_i=z_j$ or $z_i=0.$ For a permutation $\sigma \in\Sigma_r$
define the linear form on $R_{A_r}$
$$Ires_{z=0}^{\sigma}f=Res_{z_{\sigma(1)}=0}Res_{z_{\sigma(2)}=0}\cdots
Res_{z_{\sigma(r)}=0}f(z_1,z_2,\ldots, z_r)=$$
$$Res_{z_1=0}Res_{z_2=0}\cdots
 Res_{z_r=0}f(z_{\sigma^{-1}(1)},z_{\sigma^{-1}(2)},\ldots, z_{\sigma^{-1}(r)}).$$

In particular for $\sigma=id$ the  linear form $f\mapsto
Ires_{z=0}f$  defined by
$$Ires_{z=0}f$$ $$=Res_{z_1=0}Res_{z_2=0}\cdots
Res_{z_r=0}f(z_1,z_2,\ldots, z_r)$$ is  called the {\em iterated
residue.}

{\bf Remark}
\begin{itemize}
\item the  linear form $f\mapsto Ires_{z=0}^{\sigma}f$ on $R_{A_r}$  induces a linear form on $S_{A_r}$, since
it vanishes on the vector space of derivatives
$\sum_{i=1}^r\partial_i R_{A_r}$.
 \item $Ires_{z=0}^{\sigma} f_w=\delta_w^{\sigma}$.
\item the $r!$ linear forms $Ires_{z=0}^{\sigma}f,\sigma \in \Sigma_r$,
 on $S_{A_r}$ are  dual to the basis $f_w$.
\end{itemize}

Iterated residues are easier  to understand, and we will see
shortly how to use them in connection to our formulas.
 Let $w\in \Sigma_r$ and $n(w)$ be the number of elements $i$ such
that $w(i)>w(i+1)$ (this is called the number of
{\em descents} of the permutation $w$ in \cite{Stanley}).
 We denote by $C_w^+\subset C({A_r}^+)$ the
simplicial cone generated by the vectors

$$\epsilon(1)
(e_{w(1)}-e_{w(2)}),\epsilon(2)
(e_{w(2)}-e_{w(3)}),\ldots,\epsilon(r-1)  (e_{w(r-1)}-e_{w(r)}),
(e_{w(r)}-e_{r+1}),$$ where $\epsilon(i)$ is 1 or -1 depending
whether $w(i)<w(i+1)$ or not.  When $w=1$, then $C_1=C(A_r^+)$.
The following lemma is  easy to see.

\begin{lemma}\label{specialperm}
Let  $a=\sum_{j=1}^{r+1} a_je_j$ in $E_r$. The cone $C_w^+\subset
E_r$ is  given by the following system of inequalities
$\sum_{j=1}^i  a_{w(j)}\geq 0$, for all $i $ such that $w(i)<
w(i+1)$, but $\sum_{j=1}^i a_{w(j)}\leq 0$ if $w(i)>w(i+1)$.
\end{lemma}

From   Theorem \ref{formule} we obtain:

\begin{theorem}[\cite{baldonivergne}]
 Let $\c$ be a  chamber of $C(\Phi)$. Consider the set of elements $w\in \Sigma_r$
such that $\c\subset C_w^+$. Then, for $f\in S_{A_r}$,

$$\langle\langle \c,f\rangle\rangle=\sum_{w\in \Sigma_r, \c\subset
C_w^+}(-1)^{n(w)}  Ires_{z=0} w^{-1} f.$$ In particular for
$f=J_{\Phi}(a)$ we obtain

{\bf Formula 1:} for $ a\in \overline{\c}$,  we have
$$v(\Phi,a)=\langle\langle \c,J_{\Phi}(a)\rangle\rangle=\sum_{w\in \Sigma_r, \c\subset
C_w^+}(-1)^{n(w)}  Ires_{z=0}^w J_{\Phi}(a).$$
\end{theorem}

The formula is a direct consequence of the fact that $Ires_{z=0}^w$ is
the dual basis of $f_w.$ We have seen that to compute the number of
integral points of our polytope we need to compute $K_{\Phi}(a).$ Let
$t_j=m_{j,j+1}+\cdots+m_{j, r+1}-1,$ where we recall that $m_{i,j}$ is
the multiplicity of the root $e_i -e_j$ in $\Phi.$ After a change of
variable for the total residue, we obtain:

\begin{theorem}\label{functionkostant}

Let $a=\sum_{i=1}^{r+1}a_i e_i$ in $E_r\cap \Z^{r+1}$. Let
$$f_{\Phi}(a)(z)=\frac{(1+z_1)^{a_1+t_1}(1+z_2)^{a_2+t_2}\cdots
(1+z_r)^{a_r+t_r}}{z_1^{m_{1,r+1}}z_2^{m_{2,r+2}}\cdots
z_r^{m_{r,r+1}} \prod_{1\leq i< j\leq r}(z_i-z_j)^{m_{ij}}}.$$

Then
{\bf Formula 2:} for $ a\in \overline{\c}$,

$$k(\Phi,a)=\sum_{w\in \Sigma_r, \c\subset C_w^+}(-1)^{n(w)}  Ires_{z=0}^w
f_{\Phi}(a).$$

\end{theorem}

We now want to give an even more explicit formulation of the above result
  suited to be directly implemented.
  For this purpose we need to introduce some more notations. For $a\in
E_r$, let $def(a)$ be defined by   $def(a)=a+ \epsilon
\sum_{\alpha \in \Phi} \alpha +\epsilon^2( \sum_{i=1}^{r}
e_i-re_{r+1})$ with $\epsilon=\frac{1}{2mr^2}$ and $m$  the
maximum of the multiplicities $m_{ij}$.

 A {\em wall} of $A^+_r$ is a hyperplane generated by $r-1$ linearly
independent elements of $A^+_r$. The cells in $C(A^+_r) \setminus
{\cal H}$ ( ${\cal H}$ being the set of hyperplanes for $A^+_r$ ) are
open cells, interior of polyhedral cones. We will call these open
cells {\em topes}.  We will say that $a\in C(A^+_r)$ is {\em regular}
if $a$ is not on any wall for $A^+_r$. The walls of $A^+_r$ are easily
characterized since they are the kernel of a linear form as $
\sum_{i\in J} a_i $ where $J$ is a subset of $\{1,2,\ldots,r\}$.  It
is then easy to decide whether a vector $a$ is regular or not.

If $a$ is a regular element we let $\c$ denote the unique chamber of
 $C(A^+_r)$ containing it. Then the set $Sp(a)=\{w\in \Sigma_r |
 \c\subset C_w^+\}$ can be computed without explicit knowledge of the
 chamber. In fact one can easily see that the set $Sp(a)$ consists of
 those $w\in \Sigma_r$ that\ satisfy the following conditions:

\[
\left \{ \begin{tabular}{ll}
\mbox{if}
$ a_{w(1)}\geq 0$   \mbox{then}  &$w(1)<w(2)$ \ \mbox{else} $ w(1)>w(2)$  \\
\mbox{if}  $a_{w(1)}+a_{w(2)}\geq $0   \mbox{then}& $w(2)<w(3)$ \ \mbox{else} $w(2)>w(3)$\\
$\cdots$&\\
\mbox{if} \ \ $a_{w(1)}+\cdots+a_{w(i)}\geq 0$ \mbox{then}  & $w(i)<w(i+1)$ \ \mbox{else} $w(i)>w(i+1)$
\\
$\cdots$&\\
\mbox{if} \ \ $a_{w(1)}+\cdots+a_{w(r-1)}\geq 0$   \mbox{then}  &$w(r-1)<w(r)$  \mbox{else} \
$w(r-1)>w(r) $\\
\end{tabular} \right \}
\]
An element of $Sp(a)$ will be called a {\em special permutation}.

Remark that if $a_i\geq
0$ for all $i\leq r$, then $a=\sum_{i=1}^ra_i
e_i-(\sum_{i=1}^ra_i)e_{r+1}$
 belongs to the closure of the {\em
nice chamber} $\c^+$  and $Sp(a)=\{id\}.$

Now we can state
Theorem~\ref{functionkostant} as follows:
\begin{theorem}
Let $\Phi\subset A^+_r$ be a system generating $E_r.$  Let $a=\sum_{i=1}^{r+1}
a_i e_i\in E_r,\ a_{r+1}=-(a_1+\cdots+a_r),\ a_i\in \Z$ and assume that $ a\in C(A^+_r).$

 Write $$f_{\Phi}(a_1,a_2,\ldots,
a_r)(z)=\frac{(1+z_1)^{a_1+t_1}(1+z_2)^{a_2+t_2}\cdots
(1+z_r)^{a_r+t_r}}{z_1^{m_{1,r+1}}z_2^{m_{2,r+2}}\cdots
z_r^{m_{r,r+1}} \prod_{1\leq i< j\leq r}(z_i-z_j)^{m_{ij}}}.$$

Then
\begin{itemize}
\item {\bf  Formula 2A:} if $a$ is regular then
$$k(\Phi,a)=\sum_{w\in Sp(a)} (-1)^{n(w)} Ires_{z=0}^w
f_{\Phi}(a).$$
\item {\bf Formula 2B:} if $a$ is not regular then
$$k(\Phi,a)=\sum_{w\in Sp(def(a))} (-1)^{n(w)} Ires_{z=0}^w
f_{\Phi}(a).$$
\end{itemize}
\end{theorem}
{\bf Remark}  Formula 2B in the theorem follows
 by observing that the chamber containing the regular
element $def(a)$ contains  $a$ in its closure.
The deformation has to be done with care
 to deal with some {\em border} cases. The following lemma,
  that we state for completeness, shows that
the deformation with $a_i$ integers is small enough  to   take
care of such cases.

\begin{lemma}\label{deformation}
Given  $a\in C(A^+_r)\cap \Z^{r+1},$ define 
$def(a):=a+ \epsilon \sum_{\alpha \in \Phi} \alpha +\epsilon^2
(\sum_{i=1}^{r} e_i-re_{r+1}),\  \epsilon=\frac{1}{2mr^2}$ where $m$
is the maximum of the multiplicities $m_{ij}$.  Then the following holds:
\begin{itemize}

\item $def(a)$ is regular, i.e. it belongs to a chamber.
\item if $\tau$ is a tope and $a \in\tau$ then $def(a)\in \tau$
\item $a\in C(A^+_r)$ if and only if $def(a) \in C(A^+_r)$ 
\item In general if $\Phi$ is a subset of $A^+_r$, $a \notin
C(\Phi)$ if and only if $def(a)\notin C(\Phi).$
\end{itemize}
\end{lemma}

For example, we obtain the following formula for the complete network
$K_{r+1}$ on $r+1$ nodes, with excess vector $a_1,a_2,\ldots, a_r,
a_{r+1}=-\sum_{i=1}^ra_i.$ In this case, the function $k(A_r^+,a)$ is
the so-called Kostant partition function and has special importance
for the representation theory of the group $GL(r+1,\C)$.
\begin{corollary}\label{completegraph}

 For $ a\in C(A^+_r)\cap \Z^{r+1},$ the Kostant partition function
is given by:
$$k(A^+_r,a)=\sum_{w\in Sp(a')}(-1)^{n(w)}  Ires_{z=0}^ w
\frac{(1+z_1)^{a_1+r-1}(1+z_2)^{a_2+r-2}\cdots (1+z_r)^{a_r}
}{z_1\cdots z_r \prod_{1\leq i< j\leq r}(z_i-z_{j})}$$
where \[a'=\left\{\begin{matrix}
a & \mbox{\rm if a is regular}\\
def(a)&  \mbox{\rm otherwise}
\end{matrix} \right.
\]

In particular, if $a_i\geq 0$ for $1\leq i\leq r$, we have
$$k(A^+_r,a)=$$
$$Res_{z_1=0}Res_{z_2=0}\cdots
Res_{z_r=0}\left(\frac{(1+z_1)^{a_1+r-1}(1+z_2)^{a_2+r-2}\cdots
(1+z_r)^{a_r} }{z_1\cdots z_r \prod_{1\leq i< j\leq
r}(z_i-z_{j})}\right).$$

\end{corollary}

Similarly we may write a formula for the transportation polytope
$T_{m,n}(d,c)$.
\begin{corollary}
 Let   $a=\sum_{i=1}^m d_i e_i-\sum_{j=1}^{n} c_j e_{m+j},$ with $d_i$ and $c_j$ non negative integers. Then  the
number of integral points in $T_{m,n}(d,c)$ is equal to
$$\sum_{w\in Sp(a')} (-1)^{n(w)} Ires_{z=0}^ w $$ $$\times
\frac{(1+z_1)^{d_1+n-1}(1+z_2)^{d_2+n-1}\cdots (1+z_m)^{d_{m}+n-1}
(1+z_{m+1})^{-c_{1}-1}\cdots (1+z_{m+n-1})^{-c_{n-1}-1}}{z_1\cdots
z_m \prod_{1\leq i\leq m\atop 1\leq j\leq n-1}(z_i-z_{m+j})}$$
where \[a'=\left\{\begin{matrix}
a & \mbox{\rm if a is regular}\\
def(a)&  \mbox{\rm otherwise}
\end{matrix} \right.
\]

\end{corollary}
\subsection {The Algorithm for Counting Integral Flows.}

Scope of this section is a brief description of the various
algorithmic procedures that were implemented with the symbolic
language {\tt Maple} and that achieve the formula for the number of
integral points described in Theorem~\ref{functionkostant}. This
software is available at \url{www.math.ucdavis.edu/~totalresidue}. 
The initial data are an $r$ by $N$ matrix $A$ whose columns are the
elements of $\Phi$ and an element $a=\{a_1,\ldots,a_r\}\in \Z^r$ that
determines the polytope.  The ingredients that we need to compute are:

\begin{enumerate}
\item The element $a'=def(a)$ obtained by deforming the initial parameter $a$.

\item The set of permutations that appear in the formula, that is 
the  set  of special permutations\label{permutations} $Sp(a')$.

\item The residues that appear in Formula 2.
\end{enumerate}

We will discuss the ingredients for each one of these steps listing
the various algorithms that are related to the part we are describing.

First of all we want to check if our vector is in $C(A^+_r)$, that is
in the cone generated by
$\{(e_1-e_2),(e_2-e_3),\ldots,(e_{(r-1)}-e_r),e_r\}$ because otherwise
the polytope is empty and there is nothing to do. To be in the cone,
$a$ must satisfy $a_1\geq 0,\,$ $a_1+a_2\geq 0,$ \dots, \,
$a_1+a_2+\cdots+a_r\geq 0.$ The procedure {\bf check-vector} verifies
whether this is true or not. In fact because of
Lemma~\ref{deformation} we may use $def(a)$ instead of $a$ and we do
this to simplify the procedures. We compute the element $def(a)$ via
the {\tt Maple} procedure {\bf def-vector}. The vector $def(a)$ is
used in all the formulas defining $Sp(a)$ instead of $a$, whether or
not $a$ is regular. This takes care of the first part.

 For finding the subset $Sp(a)$ of $\Sigma_r$, we use the procedure
{\bf special-permutations}.  We stress that using the {\tt Maple}
function $combinat[permute]$ is impractical and does not go very far
because of memory limitations. Our approach constructs recursively the
permutations subject to our conditions, thus we save much memory in
listing only those permutations.  The set $Sp(a)$ depends strongly on
the element $a$. We do not have upper bound estimates on the subset
$Sp(a)\subset \Sigma_r$, but it seems that this set is small compared
to $\Sigma_r$.  One of the worst experimental cases for the complete
graph $K_{10}$ on $10$ nodes (the case of $A_9^+$) is the case of the
vector
$a=[30201,59791,70017,41731,58270,-81016,-68993,-47000,-43001,-20000]$
where the number $Sp(a)\subset \Sigma_9$ is $9572$, certainly much
smaller that $9!$. Experiments show that the time spent to compute
this set is rather small.

 Each permutation $w\in Sp(a)$ gives rise to the simplicial cone
$C_w^+$ containing $a$, this corresponds to a vertex of the polytope
$P(A_r^+,a)$. However, clearly the cardinality of $Sp(a)$ is much
smaller that the number of vertices of the partition polytope
$P(A_r^+,a)$. For example, for $a=[a_1,a_2,\ldots,
a_{r},-(\sum_{i=1}^r a_i)]$ with $a_i>0$, we have already remarked
that the cardinality of $Sp(a)$ is $1$, as $Sp(a)$ is reduced to the
identity permutation.

Finally, for the last step we need to compute the residue. Recall that
we need to compute

$$Ires_{z=0}^w \frac{(1+z_1)^{a_1+t_1}(1+z_2)^{a_2+t_2}\cdots
(1+z_r)^{a_r+t_r}}{z_1^{m_{1,r+1}}z_2^{m_{2,r+2}}\cdots
z_r^{m_{r,r+1}} \prod_{1\leq i< j\leq r}(z_i-z_j)^{m_{ij}}}$$ 

with $w$ one of the special permutations. Let us denote by $F$ the
function appearing in the formula above. The function $F$ is a product
of a certain number of functions. This allows us to take the residues
by introducing little by little the part of the function $F$
containing the needed variable.  To make things simpler we assume that
$w$ is the identity permutation. We start by taking the residue at
$z_r=0$ of the function
$g:=\frac{(1+z_r)^{(a_r+t_r)}}{z_r^{m_{r,r+1}}\prod_{j=1}^{r-1}
(z_j-z_r)^{m_{jr}}}$ Suppose $g_r(z_1,z_2,\dots,z_{r-1})$ is the
result.  We continue by taking the residue in $z_{r-1}$ of the
function $g_r$ multiplied by all the factors of the original function
$F$ that involve the variables $z_{r-1}$ and so on. The way we compute
the residue in one variable $z$ of a function $g(z)=F(z)/z^u$, where
$F$ is analytic, is by computing the Taylor expansion of $F$ up to the
estimate we have for the order $u$ of the pole of the function $g$ and
then taking the coefficient of $1/z$.  The argument just described is
implemented via different procedures: {\bf coeex},{\bf invi},{\bf
trunc-next-function} and {\bf RRK }.
Finally, the procedure {\bf number-kostant} adds up, with a sign (the
appropriate sign is computed using {\bf segnop}), all residues coming
from the different special permutations, thus getting Formula 2.
The procedure {\bf polynomial-kostant} computes the polynomial
$a\mapsto k(\Phi,a)$ on the chamber determined by $a$.

As we pointed out we need an uniform estimate for the order of
poles appearing. The result for the order of pole is the content
of the subsection that follows and it is implemented in procedure
{\bf E}.

\subsection {Estimates for the order of poles}
 Let $G_r$ be a Laurent polynomial in the
$r$ variables $z=(z_1,z_2,\ldots,z_r)$ and let $D_r=\prod_{1\leq
i< j\leq r}(z_i-z_j).$ We have seen that we need to compute
iterated residues of the form :
$$Res_{z_1=0}Res_{z_2=0}\cdots
Res_{z_r=0}\  G_r/D_r^m.$$ The following key lemma will handle the
situations that will appear in computing the estimate we are
looking for.

\begin{lemma}
Assume that $G_r=\frac{F(z_1,\ldots,z_r)}{(z_1 z_2 \cdots
z_r)^g}H_r(\frac{1}{z_1},\ldots,\frac{1}{z_r})$
 where $F$ is analytic and $H_r$ is a
homogeneous polynomial of degree h , then
 $$Res_{z_r=0}
G_r/D_r^m$$ is a linear combination of functions of the form
$G_{r-1}/D_{r-1}^m$ with
$$G_{r-1}=\frac{F(z_1,\ldots,z_{r-1})}{(z_1 z_2\cdots  z_{r-1})^{(g+m)}}
H_{r-1}(\frac{1}{z_1},\ldots, \frac{1}{z_{r-1}})$$ where $H_{r-1}$
is a homogeneous polynomial of  degree at most $g+h-1$ and
$F(z_1,\ldots,z_{r-1})$ is analytic.

\end{lemma}

\begin{proof}
Let us prove the lemma for a monomial $H_r=z_1^{i_1}\cdots
z_{r-1}^{i_{r-1}}z_r^{i_r}$ where $i_1,i_2,\ldots,i_r$ are non-
negative integers such that $i_1+i_2+\cdots+i_r=h$. We write
$\prod_{1\leq i\leq r-1}(z_i-z_r)^m=(z_1 z_2 \cdots
z_{r-1})^m\prod_{1\leq i\leq r-1}(1-\frac{z_r}{z_i})^m$.

The Taylor expansion of $\frac{1}{\prod_{1\leq i\leq
r-1}(1-\frac{z_r}{z_i})^m}$ at $z_r=0$ is
$$\sum_{U_1,\ldots,U_{r-1}} z_1^{-|U_1|}z_2^{-|U_2|}\cdots
z_{r-1}^{-|U_{r-1}|}z_r^{|U_1|+|U_2|+\cdots+|U_{r-1}|}$$  where
$U_s=\{j_1^s,j_2^s,\ldots,j_{m}^s\}$ varies  over the $m$ tuples
of non negative integers. Write also $F(z_1,\ldots,z_r)=\sum_k
F_k(z_1,\ldots,z_{r-1})z_r^k$. Thus we obtain

$$Res_{z_r=0}\frac{G_r}{D_r^m}=$$ $$\frac{z_1^{-i_1}\cdots
z_{r-1}^{-i_{r-1}}}{(z_1 z_2 \cdots
z_{r-1})^{g+m}}\frac{1}{D_{r-1}^m} Res_{z_r=0}
\frac{F(z_1,\ldots,z_r)}{z_r^{g+i_r}\prod_{i=1}^{r-1}(1-\frac{z_r}{z_i})^m}=$$
$$(\frac{z_1^{-i_1}\cdots
z_{r-1}^{-i_{r-1}}}{(z_1 z_2 \cdots
z_{r-1})^{g+m}}\frac{1}{D_{r-1}^m}) \times
$$$$\sum_{k=0}^{g-1+i_r}
(F_k(z_1,\ldots,z_{r-1})\sum_{U_1,\ldots,U_{r-1}:|U_1|+\cdots+|U_{r-1}|=g-1+i_r-k}
z_1^{-|U_1|}\cdots z_{r-1}^{-|U_{r-1}|})$$ For $0\leq k \leq
i_r+g-1$, the monomial $$z_1^{-i_1}\cdots z_{r-1}^{-i_{r-1}}
z_1^{-|U_1|}z_2^{-|U_2|}\cdots z_{r-1}^{-|U_{r-1}|}$$ is such that
$$i_1+\cdots +i_{r-1}+|U_1|+\cdots +|U_{r-1}|=i_1+\cdots
+i_{r-1}+i_r+g-1-k\leq h+g-1$$ and we obtain the lemma.

\end{proof}

Observe that if $F=1$ then the same proof shows that $H_r$ is
homogeneous of degree precisely $h+g-1.$ Now starting from
$G_r=\frac{F(z_1,\ldots,z_r)}{(z_1\cdots z_r)^m}$ we want to
compute

$$
Res_{z_{k+1}=0}Res_{z_{k+2}=0}\cdots Res_{z_{r-1}=0}Res_{z_{r}=0}
G_r/D_r^m.$$ Applying the lemma with $h=0$, we obtain that
$$
Res_{z_{r}=0} G_r/D_r^m$$ is a linear combination of functions of
the form $\frac{G_{r-1}}{D_{r-1}^m}$ where
$$G_{r-1}=\frac{F(z_1,\ldots,z_{r-1})}{(z_1\cdots
z_{r-1})^{2m}}H(\frac{1}{z_1},\ldots,\frac{1}{z_{r-1}})$$ and $H$
is homogeneous of degree at most $m-1$, thus at the next residue
we get again a linear combination of functions of the form
$\frac{G_{r-2}}{D_{r-2}^m}$ where
$$G_{r-2}=\frac{F(z_1,\ldots,z_{r-2})}{(z_1\cdots
z_{r-2})^{3m}}H(\frac{1}{z_1},\ldots,\frac{1}{z_{r-2}})$$ with $H$
homogeneous of degree at most $2m+m-1-1=3m-2$, so finally the last
residue in $z_{k+1}=0$ leaves  a linear combination of functions
of the form $$\frac{G_{k}}{D_{k}^m}$$ with
$$G_k=\frac{F(z_1,\ldots,z_{k})}{(z_1\cdots z_{k})^{(r-k+1)m}}
H(\frac{1}{z_1},\ldots,\frac{1}{z_{k}}).$$ Here  $H$ is
homogeneous of degree at most $\frac{(r-k)(r-k+1)m}{2}-(r-k).$ In
particular, considering $H(\frac{1}{z_1},\ldots,\frac{1}{z_{k}})$
we have the estimate on poles we were looking for.
\begin{corollary}
\begin{enumerate}
\item Let  $G_r=\frac{F(z_1,\ldots,z_r)}{(z_1\cdots z_r)^m},$ with $F$ analytic.
Then the function
$$Res_{z_{k+1}=0}Res_{z_{k+2}=0}\cdots Res_{z_{r-1}=0}Res_{z_{r}=0}
G_r/D_r^m$$ has a pole in $z_k$ of order at most
$\frac{m(r-k)(r-k+1)}{2} -(r-k)$.
\item In particular  with the notation as in Theorem~\ref{functionkostant}, if
$m=maximum_{ij}{m_{ij}}$ then the   pole in $\sigma(z_k)$ of the
function

{\small
$$Res_{z_{\sigma(k+1)}=0}\cdots
Res_{z_{\sigma(r)}=0}f_{\Phi}(a_1,a_2,\ldots,
a_r)(z)=$$$$Res_{z_{\sigma(k+1)}=0}\cdots
Res_{z_{\sigma(r)}=0}\frac{(1+z_1)^{a_1+t_1}(1+z_2)^{a_2+t_2}\cdots
(1+z_r)^{a_r+t_r}}{z_1^{m_{1,r+1}}z_2^{m_{2,r+2}}\cdots
z_r^{m_{r,r+1}} \prod_{1\leq i< j\leq r}(z_i-z_j)^{m_{ij}}}$$ }

has at most order $\frac{m(r-k)(r-k+1)}{2} -(r-k)$ independently
from $\sigma \in \Sigma_r.$

\end{enumerate}
\end{corollary}

\section{The Chamber Complex} \label{chambers}

In this section we discuss the chambers and how to compute them.  It
is important to emphasize that everything that we present in this
section is valid for general matrices, not necessarily unimodular.
There is an implementation of these ideas in the {\tt Maple} program
{\bf chambers} available at \url{www.math.ucdavis.edu/~totalresidue}.
Let $\Delta^+$ the set of distinct vectors $\{\Phi\}$. Recall the
{\em chamber complex} is the polyhedral subdivision of the cone
$C(\Delta^+)$ of nonnegative linear combinations of $\Delta^+$.
It is defined as the common refinement of the simplicial cones
$C(\sigma)$ running over all possible basic subsets $\sigma$ of
$\Delta^+.$ To be more precise we introduce now notation and the
key definitions.  In what follows, when we consider a subset
$I=\{s_1,s_2,\ldots, s_k\}$, where the elements $s_i$ of $I$ are
subsets of a set $X$, we assume there is a partial order on $I$ by
containment. Thus the set of minimal elements of $I$ is denoted by
$minimalize(I)$. We adopt the convention that the intersection of
an empty family of subsets of $X$ is $X$ itself.

Let $\Delta^+$ be the set $\{\phi_1,\phi_2,\ldots,\phi_N\}$ of
vectors in $\R^r$. Recall that a {\em wall} is a hyperplane in
$\R^r$ spanned by $(r-1)$ vectors of $\Delta^+$. Each wall $W$
partitions the set of indices $\{1,2,\ldots,N\}$ into three sets:
$zeros(W)=\{i | \phi_i \in W\}$, and two disjoint subsets
$pos(W)$, $neg(W)$ whose union $pos(W)\cup neg(W)$ is precisely
the subset of $\{1,2,\ldots,N\} \setminus zeros(W)$. We consider
the set $\{pos(W),neg(W)\}=\{neg(W),pos(W)\}$.  We denote by $\CB$
the set of  subsets $\sigma$ of $\{1,2,\ldots,N\}$ such that
$\sigma$ is of cardinality $r$ and the set of vectors $\{\phi_i | i\in
\sigma\}$ are linearly independent. For convenience, we continue
to call such  $\sigma$ a {\it basic subset} of $\Delta^+$,
thinking  of $\sigma$ as a subset of integers  or  as a subset of
elements of $\Phi$ labeled by indices.

For $\sigma\in \CB$, we consider the closed cone $C(\sigma)$ generated
by $\sigma$. If $I$ is a subset of $\CB$, let $F(I)=\cap_{\sigma\in
I}C(\sigma)$ be the intersection of the cones $C(\sigma)$, when
$\sigma$ runs in $I$. We will say that $I$ is a {\em feasible subset}
of $\CB$ if the interior of $F(I)$ is non empty. A {\em combinatorial
chamber} $I$ is a maximal feasible subset of $\CB$.  The polyhedral
cone $F(I)$ will be called a {\em geometric chamber}. The actual
chamber $Chamber(I)$ is the interior of $F(I)$. Reciprocally, the
collection $I$ is entirely determined by $F(I)$. We have
$I=\{\sigma\in \CB | F( I)\subset C(\sigma)\}.$ The collection of all
geometric chambers and their faces forms a polyhedral complex that
partitions the cone $C(\Delta^+)$, the so called {\em chamber complex}
\cite{tatiana1,billeraetal,deloeraetal}.

\begin{figure}
\begin{center}
    \includegraphics[scale=.5]{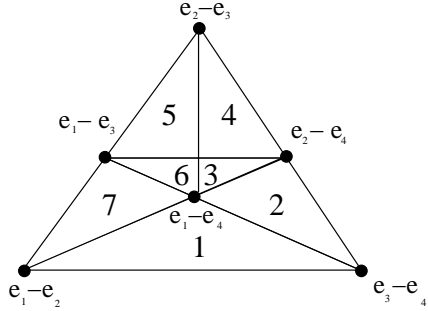}
\caption{A slice of the chamber complex for $K_4$}
\label{A4chambers}
\end{center}
\end{figure}

Figure \ref{A4chambers} shows an example, the chamber complex for
the cone associated to the acyclic complete graph $K_4$ we
discussed in the previous section.  The picture represents a
2-dimensional slice of the cone decomposition (the cone is
3-dimensional and pointed at the origin). The $6$ dots labeled
$(e_i-e_j)$ on the drawing are the intersections of the rays
$\R^+(e_i-e_j)$ with the hyperplane $(3 x_1+x_2-x_3-3x_4)=2$.
Seven chambers, numbered from 1 to 7, are present. In the
configuration of vectors of Figure \ref{A4chambers} there are
seven walls, one for each of the distinct lines obtained from the
vectors in the configuration.

\begin{figure}[htpb]
\begin{center}
    \includegraphics[scale=.5]{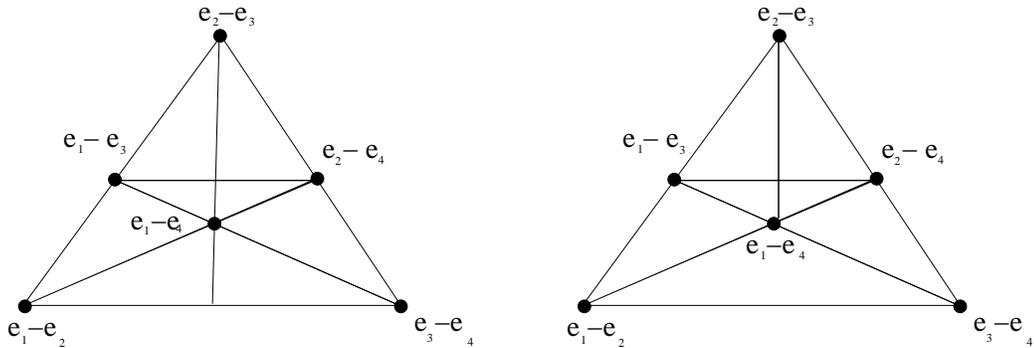}
\caption{8 topes (left) versus 7 chambers (right)}
\label{chambervshyper}
\end{center}
\end{figure}

Let ${\cal H}$ denote the hyperplane arrangement consisting of all
walls. ${\cal H}$ contains as a subset the walls of the chambers.
The cells in $C(\Delta^+) \setminus {\cal H}$ are open cells,
interior of polyhedral cones. We will call these open cells {\em
topes} (following the oriented matroid terminology
\cite{bjorneretal}). Note that the set of topes is
(typically) a much finer subdivision of $C(\Delta^+)$ than its
chambers. See Figure \ref{chambervshyper} for a comparison between
the chamber complex and the tope complex of the hyperplane
arrangement ${\cal H}$ associated with the example in Figure
\ref{A4chambers}.

A tope $\tau$ of $C(\Delta^+)$ does not touch any wall of
$\Delta^+$. Then, for each wall $W$, we denote by $pos(W,\tau)$
the set of elements $i\in \{1,2,\ldots, N\}$ such that $\phi_i\in
\Delta^+$ lies on the same open half-space determined by $W$ than
the tope $\tau$. We say that $pos(W,\tau)$ is a {\em non-face}
(this terminology is justified because these are the non-faces of
a certain simplicial complex in the sense of Chapter two of
\cite{Stanley2}). We denote by $Chamber(\tau)$ the chamber
containing the tope $\tau$.

To each tope $\tau$, we associate the family of positive non-faces
determined by the tope $\tau$ (we have a non-face for each wall).
Let us call this full family $Polarized(\tau)$. Consider the
family $MNF(\tau)$ of minimal elements of $Polarized(\tau)$. This
is the family $MNF(\tau)=minimalize(Polarized(\tau))$.  The first
main observation is that we can reconstruct the chamber
$Chamber(\tau)$ containing the tope $\tau$ from the set
$MNF(\tau)$. This is very useful to construct one initial chamber.
Later all others will be found from it.

The set $MNF(\tau)$ is a set of non-faces.  Let $f$ be the
cardinality of the set $MNF(\tau)$. Let us list
$MNF(\tau):=\{p_1,p_2,\ldots,p_f\}$. Each $p_i$ is a non-face. We
construct the family $\CP(\tau)$ of sets $\nu$ of the form
$\nu:=\{i_1,i_2,\ldots, i_f\}$ with $i_1\in p_1$, $i_2\in
p_2$,\ldots, $i_f\in p_f$. These we call {\em transversals} of a
family of sets. This family is denoted by $transversal(MNF(\tau))$
in the computer program we present. Again $\CP(\tau)$ is a set
whose elements are sets of indexes, its elements being subsets of
$\{1,2,\ldots, N\}$. The cardinality of a set $\nu\in \CP(\tau)$
may be smaller than $f$, as the family $MNF(\tau)$ does not
consists of disjoints sets. It is important to observe that if
$\nu$ is in $\CP(\tau)$, then for any wall $W$, the intersection
$\nu\cap pos(W,\tau)$ is not empty. We have the theorem.

\begin{theorem}\label{jesusrecipe}
The minimal elements of the family
$\CP(\tau):=transversal(MNF(\tau))$ are exactly the basic subsets
$\sigma$ of $\Delta^+$ such that $\tau\subset C(\sigma)$.
\end{theorem}

In other words, given the set $MNF(\tau)$ associated to a tope
$\tau$, the family of basic subsets $\sigma$ of $\Delta^+$ such
that $\tau$ is contained in $C(\sigma)$ is precisely the set
$minimalize(transversal(MNF(\tau)))$. We are going to prove this
theorem. We start by a lemma.
\begin{lemma} \label{genset}
Every $\nu\in \CP(\tau)$ is such that the set of vectors
$\{\phi_i | i\in \nu\}$  generates $\R^r$.
\end{lemma}

\begin{proof}
Let us see that a set $\nu\in \CP(\tau)$ generates $\R^r$. Indeed, if
not, the set of vectors $\{\phi_i| i\in\nu\}$ would be contained in a
wall $W$. Consider the set $pos(W,\tau)$ and a minimal element $p$ of
the family $MNF(\tau):=minimalize(Polarized(\tau))$ contained in
$pos(W,\tau)$. Then $p$ (meaning the set of elements $\phi_i$ indexed
by $p$) is contained in one of the open half-space determined by $W$.
Thus, contrary to our hypothesis, we would have $\nu\cap
p=\emptyset$. 
{\bf QED}
\end{proof}

We go on proving Theorem \ref{jesusrecipe}.

\begin{proof}
 Let $\sigma$ be a basic subset of $\Delta^+$ ($\sigma$ (elements
 indexed by $\sigma$) generates a simplicial cone). We now prove that
 if $\tau\subset C(\sigma)$, then $\sigma\in \CP(\tau)$ and is a
 minimal element in the family of tranversal sets $\CP(\tau)$.

 For each wall
 $W$, the set $\sigma\cap pos(W,\tau)$ is non empty. Otherwise
 $\sigma$ would be contained in the closed half space determined by
 $W$, but would be on the opposite to $\tau$ with respect to $W$, and
 the cone $C(\sigma)$ will not contain $\tau$. Let us pick for
 each $p\in MNF(\tau)$ an element $\phi_p\in \sigma\cap p$. It follows
 that $\sigma$ contains necessarily the set $\nu:=\{ \phi_p | \phi_p\in
 \sigma\cap p; p\in MNF(\tau)\}$, belonging to the family
 $\CP(\tau)$. But then $\sigma=\nu$, as $\sigma$ is a basic subset of
 $\Delta^+$ and $\nu$ indexes a set of generators of $\R^r$ by Lemma
 \ref{genset}. Furthermore $\sigma$ is minimal, as all sets belonging
 to the family $\CP(\tau)$ have cardinality at least equal to $r$.

We now prove the converse.  Let $\nu$ be a minimal set of
$\CP(\tau)$. We claim that $\tau$ is contained in the cone
$C(\nu)$. Otherwise, there would be a wall $W$ separating $\tau$
and $C(\nu)$. But by construction of $\nu$ there is an element
$p\in \nu$ contained in $pos(W,\tau)$; a contradiction with $W$
separating $C(\nu)$ and $\tau$. Now all we have to prove is that
$\nu$ has cardinality $r$.

Let $x$ be a point in $\tau$. By Caratheodory theorem, there is a
basic subset $\sigma$ contained in $\nu$ such that $x \in
C(\sigma)$. Then the tope $\tau$ is entirely contained in
$C(\sigma)$ because a tope is, by definition, not separated in two by
any hyperplane. The set $\sigma$ belongs to $\CP(\tau)$ by the
preceding discussion. But $ \sigma\subset \nu$ and $\nu$ is
minimal, thus $\nu=\sigma$.

So we conclude that the set $Chamber(\tau)$ of basic subsets
$\sigma$ of $\Delta^+$ such that $\tau\subset C(\sigma)$ is the
set $minimalize(\CP(\tau))$ of minimal elements of
$\CP(\tau)=transversal(MNF(\tau))$. {\bf QED}

\end{proof}

The {\em lexicographic tope} is the tope containing the vector
$\xi=\phi_1+\epsilon \phi_2+\epsilon^2 \phi_3+\cdots$ where
$\epsilon$ is a small number. The lexicographic chamber is the
chamber that contains the lexicographic tope.

\begin{corollary} \label{algorithm}
The following algorithm determines the $r$-simplicial cones
$C(\sigma)$ that contain the lexicographic chamber associated with
a particular labeling of the elements of $\Delta^+$, by finding
the basic sets $\sigma$ that define them.

\begin{enumerate}
\item Create the list $L$ of lexicographic nonfaces $pos(W,\tau)$
where $\tau$ is the lexicographic tope, and $W$ runs over all
possible walls of $\Delta^+$.

\item Let $F=\{A_1,A_2,\dots,A_m\}$ be the minimal non-faces from $L$.
\item Find the transversal sets to the family $F$ then minimalize
the set of transversals. The result is $\sigma_1,\dots,\sigma_k$
the desired basic sets.
\end{enumerate}
\end{corollary}

Now we are concerned with producing all other chambers from one
initial chamber, such as the lexicographic chamber. For this we need
to understand the polyhedron $F(I)$. This is a pointed polyhedral
cone. We recall, say from Chapter 8 in the book \cite{schrijver}, that
for a polyhedron $P$ (e.g. $F(I)$) given by a finite set of
inequalities $Ax \leq b$, a {\em supporting hyperplane} is an affine
hyperplane $\{ x | cx=d \}$ such that $d=max\{cx | Ax \leq b\}$. A
subset of $P$ is a {\em face} if $F=P$ or $F$ is the intersection of
$P$ with a supporting hyperplane of $P$.  A {\em facet} of $P$ is a
maximal face distinct from $P$.  We say a wall $W$ is an {\em
essential wall} of the geometric chamber $F( I)$, if $F(I)\cap W$ is a
facet of the pointed polyhedral cone $F(I)$. This is equivalent to $W$
being a supporting hyperplane of $F(I)$ and $dim( F( I) \cap W)=r-1.$
We say that two geometric chambers $F(I)$ and $F(I')$ are {\em
$W$-adjacent} if they share a common essential wall $W$ and $dim(F(I)
\cap F(I') \cap W)=r-1$. In particular, the wall $W$ is an interior
wall. In what follows, unless is necessary to avoid ambiguity, we will
simply refer to ``adjacent chambers'' without specifying the wall they
share. We present now an operation that allows us to move, under
certain conditions, from a geometric chamber to another adjacent
geometric chamber.  Since the geometric chambers form a connected
polyhedral complex, we can then apply some standard search procedure,
such as depth-first search, to enumerate and list all chambers.

We denote by $\CW$ the set of subsets $\nu$ of $\{1,2,\ldots,N\}$
such that $\nu$ is of cardinality $r-1$ and the set of vectors
$\{\phi_i | i\in \nu\}$ are linearly independent. In other words, if
$\nu$ is in $\CW$, the vector space $\CL(\nu)$ spanned by the
vectors $\{\phi_i | i\in\nu\}$ is a wall $W$. If $W$ is a wall we
denote by $\CW(W)$ the subset of $\CW$ with elements those $\nu$
such that $\CL(\nu)=W$.

If $\nu$ is in $\CW$, we consider the subsets $zeros(\nu)$,
$pos(\nu)$ and $neg(\nu)$. If $i$ is not  in $zeros(\nu)$ , then
$\nu \cup \{i\}$ is an element of $\CB$. We denote by
$\delta^+(\nu)$ the subset of $\CB$ consisting of elements
$\sigma=\nu \cup \{i\}$ where $i$ runs in $pos(\nu)$; denote
$\delta^-(\nu)$ the subset of $\CB$ consisting of elements
$\sigma=\nu \cup \{i\}$ where $i$ runs in $neg(\nu)$;

If $W$ is a wall, and $\sigma$ a subset of $\{1,2,\ldots,N\}$ we
denote by $\sigma \cap W=\sigma\cap zeros(W)$.
 We denote by $\CB(W|facet)$ the subset of $\CB$ consisting of
those elements $\sigma$ such that $\sigma \cap W$ is of
cardinality $(r-1)$.  In other words, $W$ is spanned by a facet of
the cone $C(\sigma)$.  We denote by $\CB(W|cut)$ the subset of
$\CB$ consisting of elements $\sigma$ such that both sets $\sigma
\cap pos(W)$ and $\sigma\cap neg(W)$ are non empty. For any subset
$I$ of $\CB$, we denote by $I(W|facet)=I\cap \CB(W|facet)$ and by
$I(W|cut)=I\cap \CB(W|cut)$.

 Let $I$ be a combinatorial chamber
which is a maximal feasible subset of $\CB$ . Let $W$ be a wall,
we define $B(W,I)=\{\sigma\cap W| \sigma\in I(W|facet)\}$. This is
a subset of $\CW(W)=\{\nu\in \CW | \CL(\nu)=W\}$. If $W$ is an
essential wall of $F(I)$, then (as we will see later)  for each
subset $\nu\in B(W,I)$ either $\delta^+(\nu)$ is contained in $I$
or $\delta^-(\nu)$ is contained in $I$, but not both. 

\begin{figure}[htpb]
\begin{center}
\includegraphics[scale=.25]{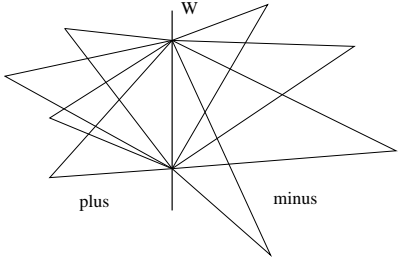}
\caption{A reflexion exchanges the simplicial cones supported
on opposite sides of a wall.}
\label{badreflex}
\end{center}
\end{figure}

If $W$ is an interior wall then define the {\em reflexion operation},
this is a new combinatorial chamber denoted by $reflexion(I,W)$. We
keep in $reflexion(I,W)$ all elements $\sigma \in I(W|cut)$, while we
replace each subset $\delta^+(\nu)\subset I(W|facet)$ by its opposite
$\delta^-(\nu)$.  The operation of reflexion has also received the
name of {\em flip} by several authors. Applying a reflexion over any
wall may not yield an adjacent chamber, as we see in the example of
Figure \ref{badreflex}

\begin{figure}[htpb]
\begin{center}
\includegraphics[scale=.3]{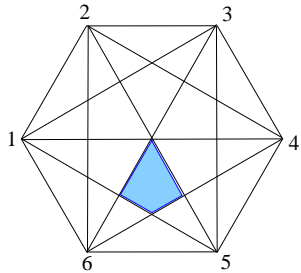}
\caption{A reflexion using the wall 1,4 does not give a chamber}
\label{badreflex}
\end{center}
\end{figure}

The important fact is that if one performs the reflexions over
essential walls the result is the desired one:

\begin{lemma}
If $W$ is an essential interior wall of $F(I)$, and let
$reflexion(I,W)$ the geometric chamber obtained by reflexion of $I$
along the essential wall $W$ . Then the set $reflexion(I,W)$ is the
combinatorial chamber associated to the $W$-adjacent chamber sharing $W$
with $F(I)$.
\end{lemma}

Clearly all elements $\sigma \in I(W|cut)$ and elements in
$\delta^-(\nu)$, when $\nu$ runs over $B(W,I)$, give rise to
simplicial cones  containing the $W$-adjacent chamber. Conversely, any
$\sigma$ in $\CB$ such that the cone $C(\sigma)$ contains the
$W$-adjacent chamber is either in $I(W|cut)$ or in a set of the form
$\delta^-(\nu)$, with $\nu\in B(W,I)$.

The above lemma stresses the importance of determining the
essential walls and that is what we describe next.  Each essential
wall $W$  is described by a linear inequality, that reaches
equality at $F(I)\cap W$. The chamber is contained in the
corresponding half-space.  The presentation we have of the chamber
is as the intersection of simplicial cones, their facets provide
us with a system of inequalities whose solution is precisely the
chamber. The trouble is that this system contains {\em redundant}
inequalities. An inequality is redundant if it is implied by the
other constraints in the system, so redundant inequalities can be
removed.

Our algorithm for finding the essential walls is based in the
following statement, which is essentially Theorem 8.1 in page 101
of \cite{schrijver}. Here we state it for full-dimensional
polyhedra (thus no equality constraints are present):

\begin{theorem} If no inequality in the system $Ax \leq b$ defining
the full-dimensional polyhedron $P$ is redundant, then there
exists a one-to-one correspondence between the facets of a
polyhedron and the inequalities in $Ax \leq b$ given by $F=\{ x\in
P | a_ix=\beta_i \}$, for any facet $F$ of $P$ and any inequality
$a_ix \leq \beta_i$ from the system $Ax \leq b$.
\end{theorem}

So if we manage to remove redundant inequalities from the original
system of inequalities associated to $F(I)$ we would have found
the essential facets of the pointed polyhedral cone $F(I)$. To do
this let us describe a direct method. Let $ A x \le b, s^T x \le
t$ be a given system of $ m+1$-inequalities in $ d$-variables $
x=(x_1,x_2,\ldots,x_d)^T$. We want to test whether the subsystem
of first $ m$ inequalities $ A x \le b$ implies the last
inequality $ s^T x \le t$. If so, the inequality $s^T x \le t$ is
redundant and can be removed from the system. A linear programming
formulation of this is rather simple:

\begin{displaymath}\begin{array}{lll} f^* = &\text{maximize} & s^T x\\ &\text{subject to} & A x \le b\\ & & s^T x \le t+1. \end{array}\end{displaymath}

Then the inequality $ s^T x \le t$ is redundant if and only if the
optimal value $ f^*$ is less than or equal to $t$.  By
successively solving this LP for each untested inequality against
the remaining inequalities, one would finally obtain an equivalent
non-redundant system. Thus the algorithm to recover all the essential
walls as follows.

\begin{enumerate}

\item Find the inequalities of each of the simplicial cones in
$F(I)$. 

\item Remove redundant inequalities using linear programming until
there is no redundant inequality left. By the previous theorem the
wall is uniquely determined by setting to equality the
inequalities.
\end{enumerate}

Thus to find all the chambers, we have

\begin{corollary} The following algorithm finds all the chambers
of the vector set $\Delta_+$:

\begin{enumerate}
\item Find the lexicographic chamber $I_{initial}$. Put that
as the first element of a list of chambers $L$.

\item Pick an element $I$ of $L$ for which we have not yet found
its adjacent chambers. Determine its essential walls $W$ using the method above.

\item Perform the reflexions $reflexion(I,W)=I(W)$ for each essential
interior wall $W$.

\item Add the $I(W)$ to the list $L$ of existing chambers if not
already there, and continue until we have found adjacent chambers for
all elements in $L$.

\end{enumerate}
\end{corollary}

Although we have a concrete algorithm now to generate all chambers
for practical reasons it is highly desirable to improve the speed
on recognizing the essential walls. For this we  prove some
necessary conditions of the essential walls of a chamber:

\begin{proposition} \label{essentialwalls}
Let $I$ be a combinatorial chamber (a maximal feasible subset of
$\CB$). Let $W$ be a wall of $\Delta^+$. If $W$ is an essential
wall of $F(I)$, then the following conditions hold true:

\begin{enumerate}
\item  $I=I(W|facet)\cup I(W|cut)$.
\item $I(W|facet)\neq \emptyset$.
\item  For each $\nu\in \CW$, either

$\delta^{+}(\nu)\cap I\neq \emptyset$.  Then
$\delta^{+}(\nu)\subset I$ and $\delta^-(\nu)\cap I=\emptyset$; or
$\delta^{-}(\nu)\cap I\neq \emptyset$. Then
$\delta^{-}(\nu)\subset I$ and $\delta^+(\nu)\cap I=\emptyset$;

\item Assume $I(W|cut)$ is not empty . Then  $\cap_{\sigma\in
I(W|cut)}\stackrel{\circ}{C(\sigma)}$ intersects $W$ in an $(r-1)$
dimensional set.
\end{enumerate}
\end{proposition}

We start the proof.  Let $\{v_1,\ldots,v_{r-1}\}$ be independent
 vectors in $\R^r$, generating a cone contained in $ F( I) \cap W$. If
 $\sigma=\{\phi_1,\ldots,\phi_r\} \in I$ we denote by $A_{\sigma}$ the
 matrix expressing $\{v_1,\ldots,v_{r-1}\}$ in terms of $\sigma$, that
 is $v_i=\sum_{j=1}^ra_{ji} \phi_j$. The matrix $a_{ji}$ has non
 negative entries for any $\sigma\in I$.

 Denote by
$$A_{i,\sigma}=[a_{i,1},\ldots,a_{i,r-1}]$$
the components of $v_1,v_2,\ldots, v_{r-1}$ on $\phi_i$. These are
the columns vectors of $A_\sigma$.

\begin{lemma}
Assume $W$ is an essential wall of $F(I)$. Suppose $W$ is spanned
by the vector set $\{v_1,\ldots,v_{r-1}\}.$ Then for each
$\sigma\in I$, either
\begin{enumerate}
\item[a)] $ A_{i,\sigma}\neq 0$ for all  $i$,

or

\item[b)] there exists an index $k$ such that $A_{k,\sigma}=0$ while
$A_{s,\sigma}\neq 0, s\neq k$.
\end{enumerate}

If $\sigma$ verifies the condition a), then $\sigma\in I(W|cut)$.
If $\sigma$ verifies the condition b), then $\sigma\in
I(W|facet)$.

\end{lemma}
\begin{proof}

Indeed, suppose that by rearranging the indices
$A_{1\sigma}=0,\ldots,A_{q\sigma}=0$, then the vectors
$\{v_1,\ldots,v_{r-1}\}$ belong to the linear space spanned by
$\{\phi_{q+1},\ldots,\phi_{r}\}$ forcing $r-q\geq r-1$ that is
$q=1$. Thus if $\sigma\in I$, $\sigma$ verifies either a) or b).

Suppose we are in the first case. We now prove that $\sigma\in
I(W|cut)$. Let us see that if all the $A_{i,\sigma}$ are non zero
vectors , then there $exists$ an element $X\in W \cap
\stackrel{\circ}{{C(\sigma)}}$ which would force that $\sigma \in
I(W|cut)$.  Let $X=\sum_{i=1}^{r-1}f_i v_i$ with $f_i>0$, then
$X\in W$. On the other hand $X=\sum_{i=1}^{r-1}f_i\sum_{j=1}^r
a_{ji} \phi_j=\sum_{j=1}^r (\sum_{i=1}^{r-1}f_i
a_{ji})\phi_j=\sum_{j=1}^r b_j \phi_j$. Because all the vectors
$A_{i,\sigma}$ are nonzero, then $b_j>0$ for all $j$, thus $X$
belongs to the open simplicial cone spanned by $\sigma$.

Consider the second case. Suppose  for simplicity that $i_{0}=1.$
Then $\{v_1,\ldots,v_{r-1}\}$ is a subset of the linear span
${\cal{L}}\{\phi_{2},\ldots,\phi_{r}\} ,$ therefore
$W={\cal{L}}\{\phi_{2},\ldots,\phi_{r}\}$ and hence $\sigma \in
I(W|facet)$.
\end{proof}

\begin{lemma}
Let $W$ be an essential wall of $F( I)$. Then the relative
interior of the cone generated by $F( I)\cap W$ (in $W$) is
contained in $\cap_{\sigma\in
I(W|cut)}\stackrel{\circ}{C(\sigma)}\cap W$.
\end{lemma}

\begin{proof}
For $X$ in the relative interior of the cone  $F( I)\cap W$ , we
choose, as in the preceding proof,  $\{v_1,\ldots,v_{r-1}\}$
independent vectors in $\R^r$, generating  a cone contained in $
F( I) \cap W$ and such that $X=\sum_{i=1}^{r-1} f_i v_i$, with
$f_i>0$. Arguing as before, we see that $X$ is in the interior of
$C(\sigma)$ for any $\sigma\in I(W|cut)$.
\end{proof}

\begin{lemma}
Let $W$ be an essential wall of $F( I)$. Then $I(W|facet)$ is not
empty.
\end{lemma}

\begin{proof}
If $I(W|facet)$ was empty, then  the open chamber $Chamber(I)$
would be $Chamber(I)=\cap_{\sigma\in
I(W|cut)}\stackrel{\circ}{C(\sigma)}$ and would intersect $W$ due
to the preceding lemma, and would not be contained on a half-space
of $W$.

\end{proof}

\bigskip

If $W$ is an essential wall of $F(I)$, then $F(I)$ is on one side
of $W$, thus we have a distinguished non face $pos(W,I)$. For each
$\nu\in \CW$ spanning $W$, we have a distinguished set
$\delta^+(\nu)=\{ \nu\cup \{i\}|\, i\in pos(W,I)\}$ of elements of
$\CB$, while $\delta^-(\nu)=\{ \nu\cup \{i\}|\, i\in neg(W,I)\}$
is disjoint from $I$:

\begin{lemma}
If $W$ is an essential wall of $F(I)$, then for every $\nu\in \CW$
such that $\CL(\nu)=W$, then

1) We have $\delta^-(\nu)\cap I=\emptyset$

2) If $\delta^+(\nu)\cap I\neq \emptyset$, then
$\delta^+(\nu)\subset I$.

\end{lemma}

\begin{proof}
Condition 1 is clear, otherwise $F(I)$ would be on the wrong side
of $W$. Now let $x\in Chamber(I)$ very closed to $W\cap F(I)$, and
$X$ in the interior of $W\cap F(I)$. Assume that $\sigma=\nu\cup
\{i\}$ belongs to $I$. Then the point $X$ is in the interior of
$C(\nu)$. The line $[x,X]$ is in the chamber $Chamber(I)$ except
at the last point $X$. It cannot cross any boundary of any
simplicial cone.  Thus we see that it stays entirely in the
interior of any simplicial cone spanned by $\nu$ and a vector
$\phi_k$ with $\phi_k$ on the same side than $x$, as clearly its
beginning $(1-\epsilon)X+\epsilon x$ is inside this simplicial
cone.
\end{proof}

Thus we have proven that if $W$ is an essential wall of $F(I)$,
the wall $W$ satisfies $1$, $2$, $3$, $4$ in the statement of
Proposition \ref{essentialwalls}.

\begin{corollary}
If $W$ is a wall of $F(I)$ satisfying  1), 2), 3) and not 4) ;
then $reflexion(I,W)$ is not a feasible subset of $\CB$.
\end{corollary}

\begin{proof}
Assume $W$ verifies 1) 2) 3).  Let $I'=reflexion(I,W)$.  If $W$
does not satisfy 4), the set $F(cut)=\cap_{\sigma\in
I(W|cut)}C(\sigma)$ does not cut $W$ in an open set. Thus is
contained in one side of the hyperplane $W$.  The set $I(W|cut)$
is left stable under the procedure reflexion.  Clearly, the other
cone $F'(facet)=\cap_{\sigma\in I'(W|facet)}C(\sigma)$ is on the
other side of the hyperplane $W$. Thus the set $I'$ is not
feasible.
\end{proof}

The following result justifies the difficulty of finding the
combinatorial chamber that contains an input vector:

\begin{proposition}
 Let $A$ be an integral matrix. Let a vector $b$ in the cone $C(A)$
generated by the columns of $A$ and a list $F$ of simplicial cones 
with rays in the columns of $A$ such that all elements of 
$F$ contain $b$. Deciding whether $F$ includes all
simplices that contain $b$, i.e. whether $F$ determines the
combinatorial chamber that contains $b$, is $NP$-hard.
\end{proposition}

\begin{proof} One well-known NP-complete problem is that of given a
complete graph with positive integral weights on the edges to decide
whether there is a hamiltonian tour of cost less than $\beta$. We will
explain now why this can be transformed of the problem of deciding
whether a list simplicial cones is already enough to determine a
chamber.

We will use a theorem by K. Murty (see Theorem 2.1 in \cite{murty}):
Consider a complete bipartite graph $K_{n,n}$. Orient the edges all in
the same direction and assign excess 1 to the tail nodes and -1 to the
head nodes of each arc. It is well known that the associated Network
polytope is the famous Birkhoff-Von Neumann polytope of doubly
stochastic matrices we saw in the introduction. This polytope is
embedded in $R^{n^2}$ and the coordinates are in correspondence with
the arcs of the bipartite network.  The associated network matrix has
rank $2n-1$, $2n$ rows and $n^2$ columns one per arc in the network
and we label them $(1,1),(1,2),\dots,(n-1,n),(n,n)$.

Extend the above network matrix by adding a row of costs, where
$c_{i,j}$, $i \not= j$,is the cost to go from $i$ to $j$, except for
the entry associated to the arc $i,i$ where one can put a huge integer
value $M$, much larger than the sum of the $n$ largest $c_{i,j}$'s. On
the righthandside of the matrix equation we add an entry of value
$\beta$. Written in terms of equations we have

$$ \sum^n_{i=1} x_{i,j}=1, \, j=1..n$$

$$ \sum^n_{j=1} -x_{i,j}=-1, \, i=1..n$$

$$ \sum^n_{i=1} \sum^n_{j=1} c_{i,j} x_{i,j}=\beta $$

$$x_{ij} \geq 0, \, for \, all \, i,j.$$

This system has now rank $2n$. The important point is: If the set of
columns $\{(1,1),(2,2),\dots, (n,n)$,
$(i_1,j_1),(i_2,j_2),\dots,(i_n,j_n)\}$ defines a simplicial cone
containing the vector $b=(1,1,1,\dots,-1,-1,-1,\beta)$, then
$(i_1,j_1),\dots,(i_n,j_n)$ must be a traveling salesman tour with
cost less or equal to $\beta$.  Thus if we take as $F$ the set of all
simplicial cones of bases that do not use all columns
$\{(1,1),(2,2),\dots, (n,n)\}$ and contain
$b=(1,1,\dots,-1,-1,\beta)$, the remaining job of deciding whether any
other cone contains the vector $b$ is then at least as hard as the
solution of the traveling salesman problem.
\end{proof}

To conclude this section it is worth mentioning that one can
abstractly apply reflexions to the non-essential walls satisfying
1) 2) and 3). The interior of the resulting ``chamber'' may
actually have empty interior in that case and thus is not useful
for us here. Nevertheless this phenomenon plays an important role
in the theory under the name of {\em virtual chambers}. In fact,
there is another characterization of the chambers using the
triangulations of the Gale diagram of the original vectors (see
\cite{ziegler} for an introduction to Gale diagrams and
triangulations).

\begin{lemma}[See \cite{billeraetal,deloeraetal}]
The face lattice of the chamber complex of a vector configuration
$A$ is anti-isomorphic to the face lattice of the secondary
polyhedron of the Gale transform of $A$, $\hat{A}$. The vertices
of the polyhedron are the regular triangulations of $\hat{A}$.
\end{lemma}

Thus generating the chambers of a network cone is the same as
generating the distinct regular triangulations of the Gale diagram
of an extended network matrix. Such calculations can be also be
done using the software {\em topcom}.

\section{Computational Experiments}

Now we present some computational experiments.  All experiments were
done in a 1 GHZ pentium computer running Linux using {\tt Maple 7}.
All our software is available at
\url{www.math.ucdavis.edu/~totalresidue}. We present our experiments
in three tables. We begin with Table \ref{k5} and Table \ref{kn} that
deal with Kostant's partition function, this is the case of acyclic
complete graphs. As we saw in Lemma \ref{reducetoacyclic}, all other
networks can be embedded into this case. We did examples in the cases
of $K_4$ $(A_3^+)$, $K_5$ $(A_4^+)$ in the first table and in the
second table we have bigger examples for the cases $A_6^+$ $A_7^+$,
$A_8^+$, $A_9^+$ and $A_{10}^+$.  We show computation times in both
tables and Table \ref{kn} also shows the cardinality of the
special permutation sets.  The computations show that the total
residue method is faster than brute force enumeration and the current
implementation of software {\tt LattE} \cite{latte} by one or two
orders of magnitud. {\tt LattE}, on the other hand, is the only
software that deals with arbitrary rational convex polyhedra.

{\tiny
\begin{sidewaystable}
\begin{tabular}{|l|r|c|} \hline
Weights on nodes  & \# of flows & secs\\ \hline
[6, 8, -5, -9] &  223 & 0.1 \\ \hline
[9, 11, -12, -8] & 330 & 0.1 \\ \hline
[1000, 1, -1000, -1] &  3002 & 0.009 \\ \hline
[4383, -886, -2777, -720] &  785528058 &  0.1  \\ \hline
[4907, 2218, -3812, -3313] & 20673947895 &  0.1 \\ \hline
[47896, 30744, -46242, -32398]  &  19470466783680  & 0.01 \\ \hline
[69295, 62008, -28678, -102625] & 179777378508547 & 0.1 \\ \hline
[3125352, 6257694, -926385, -8456661] & 34441480172695101274 & 0.01 \\ \hline
[6860556, 1727289, -934435, -7653410] & 91608082255943644656 & 0.1 \\ \hline
[12, 8, -9, -7, -4]&  14805 & 0.081 \\ \hline
[125, 50, -75, -33, -67] & 6950747024 & 0.020 \\ \hline
[763, 41, -227, -89, -488] &  222850218035543 & 0.019 \\ \hline
[11675, 88765, -25610, -64072, -10758] &  563408416219655157542748  &  0.011 \\ \hline
[78301, 24083, -22274, -19326, -60784] & 1108629405144880240444547243 &  0.029 \\ \hline
[52541, 88985, -1112, -55665, -84749] & 3997121684242603301444265332 &  0.010  \\ \hline
[71799, 80011, -86060, -39543, -26207]& 160949617742851302259767600 & 0.010  \\ \hline
[45617, 46855, -24133, -54922, -13417]  & 15711217216898158096466094 & 0.21 \\ \hline
[54915, 97874, -64165, -86807, -1817] & 102815492358112722152328 & 0.060 \\ \hline
[69295, 62008, -28678, -88725, -13900] & 65348330279808617817420057 & 0.010 \\ \hline
[8959393, 2901013, -85873, -533630, -11240903] &
6817997013081449330251623043931489475270 & 0.010 \\ \hline
[2738090,6701290, -190120, -347397, -8901863] &
277145720781272784955528774814729345461 & 0.010 \\ \hline
[6860556, 1727289, -934435, -818368, -6835042] &
710305971948234346520365668331191134724 & 0.060 \\ \hline

\end{tabular}
\centering
\caption{Testing for the complete graphs $K_4$ and $K_5$. Time is
given in seconds. Excess vectors are in the first column.} \label{k5}
\end{sidewaystable}
}

{\tiny
{\setlength{\textwidth}{40cm}
\begin{sidewaystable}
\begin{tabular}{|l|p{1.5in}|c|c|} \hline
Weights on nodes  & \# of flows & secs &  $\vert Sp(a) \vert$ \\ \hline
[1,2,3,4,5, -15 ] & 5880 & 0.02 & 1 \\ \hline
[21128,45716,79394,-76028,-31176,66462,-105496] & 58733548560911702671
\newline 16780821466940568432 \newline 553474831987566395925 & 0.22 & 8 \\ \hline
[82275,33212, 91868, -57457,47254,-64616,94854,-227390] & 22604049468113537772 \newline 228176193404009135\newline6424181 & 2.14 & 26 \\ \hline
[31994,-12275, 55541, 72295,26697,-3212,-38225,6916,-139731] & 11446847479255704222 \newline 87042245223206779226 \newline 01568734727431018393 \newline 069006356672309031382\newline51984519069399479632\newline6644137066000 & 7.94 & 24 \\ \hline
[12275,55541, 72295, 26697,-3212,-38225,6916,92409,9528, -234224] & 12970047729476531166\newline58326881685949118367\newline16319862924094634125\newline27856414458487356258\newline66474206451882923253\newline41990044115208492747\newline58896993761880000897\newline382293730 & 21.31 & 16 \\ \hline

[1,2, 3,4,5,6,7,8,9,10, -55] & 38883505145515430400 & 5 & 1 \\ \hline
[46398,36794, 92409,-16156,29524,-68385,93335,50738,75167, -54015, -285809 ] & 20889867895116832060\newline28578373441423712122\newline50684806890637191792\newline33590765780756053509\newline92237184823590262176\newline29560725791309259479\newline21077842421668832691\newline54404688022155977982\newline34585056426719876125\newline028873152 & 2193.23 & 322 \\ \hline
\end{tabular}
\centering
\caption {Testing for complete graphs $K_n$ with $n=6,7,8,9,10,11$. Time is given in seconds.} \label{kn}
\end{sidewaystable}
}
}

As it is clear on the two first tables, the computation time does not
increase significantly when the weights on nodes are very large. In
contrast, computation time becomes quickly very large, when the number
of nodes on the graph is growing. In the second table it is evident
that for a fixed number of nodes, time of computation depends strongly
of the cardinality of the set $Sp(a)$, i.e. the signs of weights on
the nodes (when all weights are positive, except the last, the
cardinality of $Sp(a)$ is $1$).

Let us stress that one of the features of our method is that it can
directly compute the polynomial $k_\Phi(a)$ giving the number of
lattice points in the polytope $P(\Phi,a)$ in the chamber determined
by $a$. In particular, the Ehrhart polynomial of the polytope
$P(\phi,a)$, i.e. the function $t\mapsto k_\Phi(t a)$ is also computed
easily from our algorithm. For example, corresponding to the first
line of Table \ref{kn}:
$$k_{A_r^+}(t,2t,3t,4t,5t,-15t)=
\frac{1}{120960}(6t+1)(t+4)(t+3)(t+2)(t+1)\times$$
$$(64921t^5+233897t^4+307649t^3+184639t^2+50574t+5040)$$
which was computed in $0.55$ seconds. In contrast,
the polynomial function $k_\Phi(a_1,a_2,a_3,a_4,a_5)$ (with $
a_5=-(a_1+a_2+a_3+a_4)$) in the chamber 
chamber $\{a_1>0, a_2>0, a_3>0, a_4>0\}$ is computed in $0.48$ seconds.

The Ehrhart polynomials for the second, third and fourth examples in
Table \ref{kn}, i.e.
$k_{A_r^+}(21128*t,45716*t,79394*t,-76028*t,-31176*t,66462*t)$,
$k_{A_r^+}(82275t,33212t, 91868t, -57457t,47254t,-64616t,94854t)$, and
$k_{A_r^+}(31994t,-12275t, 55541t, 72295t,
26697t,-3212t,-38225t,6916t)$, were computed in $1.36$ seconds,
$18.54$ seconds, and $93.36$ seconds respectively.  It is also amusing
to check the program on the value of the Kostant partition for $A_r^+$
on the vector $a=[1,2,3,4,\cdots, r, -r(r+1)/2]$. As proven by
Zeilberger \cite{zeilberger}, this value is $\prod_{i=1}^r
\frac{(2i)!}{i!(i+1)!}$.

The last table is dedicated to $4 \times 4$ transportation matrices.
In the case of transportation polytopes, i.e. complete bipartite
graphs. Here we also able to compare our speed to the special purpose
$C^{++}$ program written by Beck and Pixton \cite{beckpixton}.  Both
{\tt LattE} and Beck-Pixton's software are faster than our {\tt Maple}
implementation, with Beck-Pixton's significantly so, but it must still
be emphasized that our calculations for transportation polytopes makes
use of the fact that they are embedded inside the complete graph for
large enough number of nodes. For example the case of $4\times 4$
transportation polytopes is treated via the complete graph $K_8$.
The same kind of embedding can be done for other networks.

{\tiny
\begin{sidewaystable}
\centering
\begin{tabular}{|p{2.7in}|r|c|} \hline
Margins  & \# of lattice points & secs\\ \hline
[220, 215, 93, 64], \newline
[108, 286, 71, 127] &  1225914276768514 & 5.04 \\ \hline
[109, 127, 69, 109], \newline [119, 86, 108, 101] & 993810896945891 & 10.43 \\ \hline
[72, 67, 47, 96], \newline [70, 70, 51, 91] &  25387360604030 & 6.5 \\ \hline
[179909, 258827, 224919, 61909], \newline [190019, 90636, 276208, 168701] &  13571026063401838164668296635065899923152079 &  5.87  \\ \hline
[229623, 259723, 132135, 310952],\newline  [279858, 170568, 297181, 184826] &646911395459296645200004000804003243371154862 & 16.1  \\ \hline
[249961, 232006, 150459, 200438],\newline  [222515, 130701, 278288, 201360] & 319720249690111437887229255487847845310463475 &  16.1  \\ \hline
[140648, 296472, 130724, 309173], \newline [240223, 223149, 218763, 194882]&322773560821008856417270275950599107061263625 &  11.7 \\ \hline
[65205, 189726, 233525, 170004],\newline  [137007, 87762, 274082, 159609]  & 6977523720740024241056075121611021139576919 & 9.0\\ \hline
[251746, 282451, 184389, 194442], \newline [146933, 239421, 267665, 259009] & 861316343280649049593236132155039190682027614 & 15 \\ \hline
[138498, 166344, 187928, 186942], \newline [228834, 138788, 189477, 122613] & 63313191414342827754566531364533378588986467 & 19.4 \\ \hline
[20812723, 17301709, 21133745, 27679151], \newline [28343568, 18410455,
19751834, 20421471] &
665711555567792389878908993624629379187969880179721169068827951
& 15.6 \\ \hline
[15663004, 19519372, 14722354, 22325971], \newline [17617837, 25267522,
20146447, 9198895] &
63292704423941655080293971395348848807454253204720526472462015 &
27.4
\\
\hline

[13070380, 18156451, 13365203, 20567424], \newline [12268303, 20733257,
17743591, 14414307] &
43075357146173570492117291685601604830544643769252831337342557 &
14.8
\\
\hline

\end{tabular}
\caption{Testing for $4 \times 4$ transportation polytopes.}
\label{transport}
\end{sidewaystable}
}

If we consider the case of $4$ times $5$ matrices with weights on
nodes $[3046,5173,6116,10928], [182,778,3635,9558,11110]$, the
number of lattice points is $23196436596128897574829611531938753$
calculated in $11.15$ seconds. The number of special permutations
for this vector is $540$ while the number of vertices of the
corresponding polytope is 912. These same example takes 7.8 seconds
in {\tt LattE} and 0.1 seconds in Beck-Pixton program.

Ehrhart polynomial $k_{\Phi_{4,5}}(
(3046*t,5173*t,6116*t,10928*t,-182*t,-778*t,-3635*t,-9558*t,
-11110 *t)$ is computed in  $30.72$ seconds.

If we consider the case of $5$ times $5$ matrices with weights on
nodes
$[30201,59791,70017,41731,58270],[81016,68993,47000,43001,20000],$ the
number of lattice points is
$$24640538268151981086397018033422264050757251133401758112509495633028,
 $$
which we computed in $23$ minutes. The number of special
permutations needed is $9572$ while the number of vertices of the
corresponding polytope is $13150$. This example took 20 minutes with
{\tt LattE} and just 4 seconds with Beck-Pixton program.

Transportation polytopes were treated by Beck and Pixton
\cite{beckpixton} in a special purpose $C^{++}$ program dedicated for
this particular family of flow polytopes. Their computation is also
via residues and is the fastest at the moment. It is important to
remark that their use of residues is quite different from ours; our
main theorem can be thought of as a multidimensional analogue of the
fact that sums of the residues of a rational function on $P_1(\C)$ is
zero. It is to be expected that in a forthcoming $C^{++}$
implementation the timings discussed here will be considerable faster
than those from this preliminary {\tt Maple} implementation. Besides
obvious implementation speed ups, the ideas presented in this paper
could still be improved when the total residue method is
applied directly to the bipartite graph, not as a subnetwork of $K_n$.


\begin{thebibliography}{99}

\bibitem{tatiana1}{\bf Alekseyevskaya T.V., Gel'fand I.M., and
Zelevinsky A.} {\em Arrangements of real hyperplanes and the
associated partition function}, Soviet Math. Doklady 36, 1988,
589-593.

\bibitem{baldonivergne}{\bf Baldoni-Silva W. and Vergne M.}{\em
Residues formulae for volumes and Ehrhart polynomials of convex
polytopes.}  manuscript 81 pages, 2001. available at math.ArXiv,
CO/0103097.

\bibitem{BarviPom} {\bf Barvinok A. and Pommersheim J.}, {\em An
algorithmic theory of lattice points in polyhedra}, in: {\sl New
Perspectives in Algebraic Combinatorics} (Berkeley, CA,
1996-1997), 91-147, Math. Sci. Res.  Inst. Publ. 38, Cambridge
Univ. Press, Cambridge, 1999.

\bibitem{beckpixton}{\bf Beck M. and Pixton, D.}{\em The Ehrhart
polynomial of the Birkhoff polytope}, to appear in {\sl Discrete
and Computational Geometry}. Available at math.ArXiv, CO/0202267


\bibitem{billeraetal}{\bf Billera L.J., Gel'fand I.M, and Sturmfels
B.}{\em Duality and minors of secondary polyhedra} {\sl J. of
Comb. Theory, Ser.~B}, {\bf 57}, 1993, 258--268.

\bibitem{bincer}{\bf Schmidt J.R. and Bincer A.M.}, {\em The Kostant
partition function for simple Lie algebras.} {\sl J. Math. Phys.}
{\bf 25} (1984), no.~8, 2367--2373.


\bibitem{bjorneretal}{\bf Bj\"orner A., Las Vergnas M., Sturmfels B.,
N.~White and G.~Ziegler} {\em Oriented Matroids}, Cambridge
University Press, 1992.


\bibitem{brionvergne}{\bf Brion M. and Vergne M.} {Residue formulae,
vector partition functions and lattice points in rational
polytopes}, {\sl J. Amer. Math Soc.} 10,4 (1997), 797-833.

\bibitem{brionvergneI}{\bf Brion M. and Vergne M.} {Arrangements of
hyperplanes I: Rational functions and Jeffrey-Kirwan residue}
available at math.ArXiv, DG/9903178.

\bibitem{deloeraetal}{\bf De Loera J.A, Hosten S, Santos F., and
Sturmfels B.}, {\em The polytope of all triangulations of a point
configuration}, Doc. Math. J. DMV 1 (1996) 103--119.

\bibitem{latte}{\bf De Loera J.A., Hemmecke R., Tauzer J.,and
Yoshida R.} {\em Effective lattice point enumeration in rational
convex polytopes} available at \url{www.math.ucdavis.edu/~deloera}.

\bibitem{persi} {\bf Diaconis P. and Efron B.} {\em Testing for
independence in two-way tables: New interpretations of the
chi-square statistic}, {\sl Annals of Statistics} 13, 845-847.

\bibitem{GareyJohnson} {\bf Garey M.R. and Johnson S.J.} {\em Computers
and Intractability: A Guide to the Theory of NP-Completeness},
Freeman, San Francisco, 1979.


\bibitem{Jaeger} {\bf Jaeger F.} {\em Flows and generalized coloring
theorem in graphs} J. Combin. Theory Ser. B, 26, (1979), 205--216.

\bibitem{kirillov}{\bf Kirillov A. N.}{\em Ubiquity of Kostka
polynomials}{ available electronically at
http://front.math.ucdavis.edu/math.QA/9912094}

\bibitem{murty}{\bf Murty K.} {\em A fundamental problem in linear
inequalities with applications to the traveling salesman problem.}
{\sl Mathematical Programming} vol 2. 1972, 296-308.


\bibitem{pitmanstanley}{\bf Pitman J. and Stanley R.P.} {\em A
polytope related to empirical distributions, plane trees, parking
functions, and the associahedron} {\sl Discrete and Computational
Geometry}, 27 (2002), 603-634.


\bibitem{topcom} {\bf Rambau J.} {\em TOPCOM (triangulations of point
configurations and oriented matroids)}, software available at
http://www.zib.de/rambau/TOPCOM.html


\bibitem{schrijver} {\bf Schrijver A.} {\em Theory of Linear and
Integer Programming} Wiley series in Discrete Mathematics and
Optimization, 1982.

\bibitem{Stanley} {\bf Stanley R. P.} {\em Enumerative Combinatorics
} volume I, Cambridge Univ. Press, Cambridge, 1999.

\bibitem{Stanley2} {\bf Stanley R. P.} {\em Combinatorics and
Commutative Algebra} Birkh\"auser Boston, second edition 1996.

\bibitem{sturmfels}{\bf Sturmfels B.}{\em On vector partition
functions},{\sl J.~of Combinatorial Theory, Ser.~A} {\bf 72}
(1995) 302--309.

\bibitem{Szenesvergne}{\bf Szenes A. and Vergne M.}, {\em {Residue
formulae for vector partitions and Euler-MacLaurin sums}.}
preprint (2002), 52 pages. Available at math.ArXiv, CO/0202253.


\bibitem{zeilberger} {\bf Zeilberger, D.}, A conjecture of Chan,
Robbins, and Yuen. Available at math.ArXiv, CO/9811108 (1998).


\bibitem{ziegler} {\bf Ziegler G.} {\em Lectures on Polytopes},
Springer Verlag Graduate Texts, Berlin, 1995.


\end{thebibliography}
\end{document}